\documentclass[a4paper,11pt,reqno,twoside]{amsart}

\usepackage{amsmath}
\usepackage{amsfonts}
\usepackage{amssymb}
\usepackage{amsthm}
\usepackage{color}
\usepackage{ifpdf}
\usepackage{array}
\usepackage{url}
\usepackage{multirow,bigdelim}
\usepackage{ascmac}

\theoremstyle{oupplain}
\newtheorem{theorem}{Theorem}[section]
\newtheorem{proposition}[theorem]{Proposition}
\newtheorem{lemma}[theorem]{Lemma}
\newtheorem{corollary}[theorem]{Corollary}
\theoremstyle{oupdefinition}

\theoremstyle{oupremark}
\newtheorem{remark}[theorem]{Remark}

\numberwithin{equation}{section}

\newcommand*{\C}{\mathbb{C}}
\newcommand*{\R}{\mathbb{R}}
\newcommand*{\Z}{\mathbb{Z}}
\newcommand*{\N}{\mathbb{N}}

\newcommand{\comment}[1]{}
\title[On the Hilbert space derived from the Weil distribution]%
      {On the Hilbert space derived from \\ the Weil distribution} 
\author[M. Suzuki]{Masatoshi Suzuki}
\date{Version of \today}
\AtBeginDocument{%
\begin{abstract}
We study the Hilbert space obtained by completing 
the space of all smooth and compactly supported functions on the real line 
with respect to the hermitian form arising from the Weil distribution 
under the Riemann hypothesis. 
It turns out that this Hilbert space is isomorphic to a de Branges space 
by a composition of the Fourier transform and a simple map.
This result is applied to state new equivalence conditions for the Riemann hypothesis 
in a series of equalities.
\end{abstract}
\maketitle
}
\begin{document}

%
\section{Introduction} 
%

The Weil distribution is a distribution associated with the Riemann zeta-function $\zeta(s)$. 
Let 
\[
\xi(s) = \frac{1}{2}s(s-1)\pi^{-s/2}\Gamma\left(\frac{s}{2}\right)\zeta(s)
\]
be the Riemann xi-function, where $\Gamma(s)$ is the gamma-function. 
Let $\Gamma$ be the set of all zeros of $\xi(1/2-iz)$ without multiplicity 
and let $m_\gamma$ denote the multiplicity of $\gamma \in \Gamma$. 
The Riemann hypothesis  (RH, for short) claims that 
all nontrivial zeros of $\zeta(s)$ lie on the critical line $\Re(s)=1/2$. 
It is equivalent to the assertion that 
all $\gamma\in\Gamma$ are real.

The Weil distribution is the linear functional 
$W:C_c^\infty(\R)\to\C$ 
defined by 
\begin{equation*}
C_c^\infty(\R) \ni \psi~\mapsto~
W(\psi):=\sum_{\gamma \in \Gamma} m_\gamma \widehat{\psi}(-\gamma), 
\end{equation*}
where $C_c^\infty(\R)$ is the space of all smooth and compactly supported functions on $\R$ 
and 
\begin{equation} \label{EQ_101}
\widehat{\psi}(z) := (\mathsf{F}\psi)(z) := \int_{-\infty}^{\infty} \psi(x) \, e^{izx} \, dx
\end{equation}
is the Fourier transform. 
We omit the description of the topology of $C_c^\infty(\R)$, since we do not need it later. 
Weil~\cite{We52} (see also the note in \cite[Section 3.2]{Su22}) 
discovered that the RH is true if and only if 
the Weil distribution $W$ is nonnegative definite, that is,  
\begin{equation*}
W(\psi \ast \widetilde{\psi}) \geq  0 \quad \text{for every}~\psi \in C_c^\infty(\R),
\end{equation*}
where 
\begin{equation*} 
(\phi\ast\psi)(x):=\int_{-\infty}^{\infty} \phi(y)\psi(x-y) \, dy 
\quad \text{and} \quad 
\widetilde{\psi}(x) := \overline{\psi(-x)}. 
\end{equation*}
Further, if the RH is true, the Weil distribution is positive definite, that is, 
$W(\psi \ast \widetilde{\psi}) > 0$  for every nonzero $\psi \in C_c^\infty(\R)$. 

Using the Weil distribution, 
we define the hermitian form $\langle\cdot,\cdot\rangle_W$ on $C_c^\infty(\R)$ by 
\begin{equation} \label{EQ_102}
\langle \psi_1, \psi_2 \rangle_W 
= W(\psi_1 \ast \widetilde{\psi_2})
= \sum_{\gamma \in  \Gamma} m_\gamma 
\widehat{\psi_1}(-\gamma) (\widehat{\psi_2})^\sharp(-\gamma), 
\quad \psi_1,\,\psi_2 \in C_c^\infty(\R),
\end{equation}
where 
\[
F^\sharp(z):=\overline{F(\bar{z})}
\] 
for complex-valued functions of a complex variable. 
We often use this $\sharp$ notation. 
We call this hermitian form the {\it Weil hermitian form}. 
Yoshida \cite{Yo92} has studied the Weil hermitian form in detail 
by restricting it to a function space on a finite interval $[-a,a]$ ($a>0$). 
The subject of the present article is the behavior of 
the Weil hermitian form on the whole line $\R$. 
Yoshida proposed a method to complete a function space on a finite interval 
with respect to the Weil hermitian form without assuming the RH, 
but it does not extend to the whole line. 
\medskip

Suppose that the RH is true. 
Then the Weil hermitian form $\langle\cdot,\cdot\rangle_W$ is positive definite on $C_c^\infty(\R)$. 
Therefore, 
the completion $\mathcal{H}_W$ of the pre-Hilbert space $C_c^\infty(\R)$ 
with respect to $\langle\cdot,\cdot\rangle_W$ is defined. 
The first main result is an explicit description of the Hilbert space $\mathcal{H}_W$. 
The elements of $\mathcal{H}_W$ are equivalence classes of Cauchy sequences 
with respect to $\langle\cdot,\cdot\rangle_W$, 
where two Cauchy sequences are equivalent 
if their difference converges to zero with respect to $\langle \cdot,\cdot \rangle_{W}$. 
The representative of each class can be chosen from $L^2(\R)$ (Theorem  \ref{thm_5_1} below). 
Such a result is expected from Lemmas 2 and 3 in \cite{Yo92}. 
Therefore, we denote the class represented by $\psi \in L^2(\R)$ as $[\psi]$ 
and often identify $\psi$ with $[\psi]$. 

For the concrete description of $\mathcal{H}_W$, 
we use a de Branges space and a model space. 
The entire function $E_\xi$ defined by 
\begin{equation} \label{EQ_103}
E_\xi(z):=\xi(1/2-iz)+\xi'(1/2-iz) 
\end{equation}
belongs to the Hermite--Biehler class under the RH \cite[Theorem 1]{La06} 
and hence it defines the de Branges space $\mathcal{H}(E_\xi)$, 
where the dash on the right-hand side of \eqref{EQ_103} means differentiation 
of $\xi(s)$ with respect to $s$. 
Furthermore, the meromorphic function 
\begin{equation} \label{EQ_104}
\Theta_\xi(z):= E_\xi^\sharp(z)/E_\xi(z)
\end{equation}
in $\C$ is a meromorphic inner function in the upper-half plane $\C_+=\{z \,|\, \Im(z)>0\}$ 
under the RH, and therefore it defines the model space $\mathcal{K}(\Theta_\xi)$.  
These two Hilbert spaces $\mathcal{H}(E_\xi)$ and $\mathcal{K}(\Theta_\xi)$ 
are isomorphic with 
$\Vert E_\xi F \Vert_{\mathcal{H}(E_\xi)}
=
\Vert F \Vert_{\mathcal{K}(\Theta_\xi)}
$ 
for  every $F \in \mathcal{K}(\Theta_\xi)$
(see Section \ref{section_2} for details on 
the Hermite--Biehler class, 
de Branges spaces, and model spaces). 
Then the first result is stated as follows. 

\begin{theorem} \label{thm_1_1} 
Assume that the RH holds. 
Let $\mathcal{H}_W$, $\mathcal{H}(E_\xi)$, and $\mathcal{K}(\Theta_\xi)$ be Hilbert spaces as above. 
Then,  the map $\mathcal{K}(\Theta_\xi) \to \mathcal{H}_W$ defined by  
\[
F~\mapsto~[\psi_F], \quad \psi_F := \mathsf{F}^{-1}(F) 
\]
is an isomorphism between Hilbert spaces 
and satisfies 
\[
\Vert E_\xi F \Vert_{\mathcal{H}(E_\xi)}^2 =
\Vert F \Vert_{\mathcal{K}(\Theta_\xi)}^2 
= \pi \langle \psi_F,\psi_F \rangle_W
= \pi \langle [\psi_F], [\psi_F] \rangle_W
\]
for $F \in \mathcal{K}(\Theta_\xi)$, 
where $\mathsf{F}^{-1}$ is the Fourier inversion on $L^2(\R)$. 
\end{theorem}

This result is proved in Section \ref{section_5}. 
Note that Theorem\ref{thm_1_1} provides an isomorphism as a Hilbert space, 
not as a reproducing kernel Hilbert space. 
The space $\mathcal{H}_W$ is a space of equivalence classes of functions, 
not a space of functions.

Lagarias suggested after Theorem 1 of \cite{La06} that the norm of 
the de Branges space $\mathcal{H}(E_\xi)$ 
and the Weil hermitian form  
(the spectral side of the ``explicit formula'' of prime number theory) 
are similar. Theorem \ref{thm_1_1} shows that they are naturally coincident. 
Hence, $\mathcal{H}_W$ and $\mathcal{H}(E_\xi)$ must have 
an ``arithmetic structure'' through the  Weil explicit formula \eqref{EQ_303} below, 
but we will not discuss this further.

Connes, Consani, and Moscovici \cite[Section 4.8]{CCM24} also describe the relation 
between the theory of de Branges spaces and the Weil hermitian form, 
but their de Branges spaces $\mathcal{B}_\lambda^S$ 
and $\mathcal{H}(E_\xi)$ have completely different properties. 
Due to the difference in the generators of the de Branges spaces, 
they are not isomorphic, 
and a more obvious difference is that 
they have different spectral properties (see the second half of Section \ref{section_6}).
\medskip

One of the remarkable properties of de Branges spaces is the structure of subspaces. 
The set of all de Branges subspaces of a given de Branges space 
is totally ordered by set-theoretical inclusion (see \cite[pp. 500--506]{Wo15} for details). 
Such a structure also comes to $\mathcal{H}_W$ 
through the isomorphism of Theorem \ref{thm_1_1} 
as stated in Theorem \ref{thm_5_3} below. 

Another notable property of de Branges spaces 
is the explicit description of the family of self-adjoint extensions of 
the multiplication operator by an independent variable $F(z) \mapsto zF(z)$. 
It enables us to interpret the set of zeros $\Gamma$ 
as the set of eigenvalues of  a self-adjoint operator on $\mathcal{H}_W$. 
This means that one of the Hilbert--P{\'o}lya spaces 
is the Hilbert space $\mathcal{H}_W$ naturally obtained from the Weil distribution 
(see Sections \ref{sec_2_3} and \ref{section_6} for details). 

\medskip 

As stated in Theorem \ref{thm_1_1}, 
the Hilbert space
$\mathcal{H}_W$ is isomorphic to a de Branges space under the RH. 
Moreover, representatives of classes in $\mathcal{H}_W$ 
can be chosen from the concrete subspace $V(0)$ of $L^2(\R)$ defined in \eqref{EQ_502} below. 
It is surprising that such an explicit description of $\mathcal{H}_W$ is possible, 
and interesting in itself. 
However, it is a matter of concern that 
it is not even possible to define 
$\mathcal{H}_W$, $\mathcal{H}(E_\xi)$, and $\mathcal{K}(\Theta_\xi)$
without assuming the RH.
Fortunately, by considering a screw line of the screw function attached to $\zeta(s)$, 
which will be explained  in Sections \ref{sec_2_1} and \ref{section_4_2}, 
we can {\it unconditionally} construct two special Hilbert spaces 
$\mathcal{H}_0$ and $\mathcal{K}_0$ (in Section \ref{section_3_3}) 
to be isomorphic to 
$\mathcal{H}_W$ and $\mathcal{K}(\Theta_\xi)$, respectively, under the RH 
(Theorem \ref{thm_5_2}). 
The construction of such spaces leads to an equivalence condition for the RH stated below. 
That is the second main result. 

In Selberg's answer to the second question in \cite[p. 632]{BS08}, 
he states that the construction of a space assuming the RH will not be useful for attacking the RH. 
However, $\mathcal{H}_0$ and $\mathcal{K}_0$ 
may be useful in future research on the RH, 
since they are defined {\it without} the RH.

Let $L^2(\R)$ be the usual $L^2$-space on the real line with respect to the Lebesgue measure. 
We define 
\begin{equation} \label{EQ_105}
\mathfrak{S}_t(z) := \frac{i(1+\Theta_\xi^\sharp(z))}{2}\, \mathfrak{P}_t(z)
\end{equation}
with
\begin{equation} \label{EQ_106}
\aligned 
\mathfrak{P}_t(z)
& := \frac{4(e^{t/2}-1)}{1+2iz} + \frac{4(e^{-t/2}-1)}{1-2iz} \\ 
& \quad 
+ \frac{e^{-izt}-1}{iz}\frac{\zeta'}{\zeta}\left( \frac{1}{2}-iz \right)
+ \sum_{n \leq e^t} \frac{\Lambda(n)}{\sqrt{n}} \frac{e^{-iz(t-\log n)}-1}{iz}
\\
& \quad 
-\frac{1}{2iz} 
\left[
\frac{\Gamma'}{\Gamma}\left(\frac{1}{4}-\frac{iz}{2}\right)
-
\frac{\Gamma'}{\Gamma}\left(\frac{1}{4}\right)
\right] \\
& \quad 
- \frac{1}{2iz} 
e^{-t/2} \left[
\Phi(e^{-2t},1,\tfrac{1}{2}(\tfrac{1}{2}-iz)) 
- \Phi(e^{-2t},1,\tfrac{1}{4}) \right]
\endaligned 
\end{equation}
for a nonnegative real number $t$ and a complex number $z$, 
where $\Lambda(n)$ is the von Mangoldt function defined by 
$\Lambda(n)=\log p$ if $n=p^k$ with $k \in \Z_{>0}$ 
and $\Lambda(n)=0$ otherwise, 
and 
\[
\Phi(z,s,a) = \sum_{n=0}^{\infty} \frac{z^n}{(n+a)^s}
\]
is the Hurwitz--Lerch zeta-function. 
For negative $t$, we set $\mathfrak{S}_t(z):=\mathfrak{S}_{-t}(z)$. 
The definition of $\mathfrak{P}_t(z)$ is quite complicated.  
However, using the set $\Gamma$ of zeros of $\xi(1/2 - i z)$, 
it can be expressed in the simple form \eqref{EQ_302} 
(see Proposition \ref{prop_301} below). 
Nevertheless, as a tool for stating an equivalent condition for the RH, 
it seems preferable to have a representation that does not involve $\Gamma$. 
Thus, here we adopt a version of \eqref{EQ_302} rewritten without $\Gamma$ 
using Weil's explicit formula \eqref{EQ_303}. 
In Weil's explicit formula, the first, second, and the third-fourth lines on the right-hand side of \eqref{EQ_106} 
correspond to the poles of the completed zeta-function $\pi^{-s/2}\Gamma(s/2)\zeta(s)$, 
the non-archimedean part (Euler product), 
and the archimedean part (gamma factor), respectively. 

For this $\mathfrak{S}_t$, we first obtain the following. 

\begin{proposition} \label{prop_101}
For any fixed $t \in \R$, 
$\mathfrak{S}_t(z)$ belongs to $L^2(\R)$ as a function of $z$.
\end{proposition}
\begin{proof}
See Section \ref{sec_3_2}. 
\end{proof}

%
From this result, the mapping $t \mapsto \mathfrak{S}_t(z)$ from $\R$ to $L^2(\R)$ is defined.
By the uniformity of the $L^2$-norm of $\mathfrak{S}_t(z)$ 
on a compact set of $t$ obtained in the proof of Proposition \ref{prop_101} 
and Minkowski's integral inequality, the following holds. 

\begin{proposition} \label{prop_102} 
For $\phi \in C_c^\infty(\R)$, 
we define 
\begin{equation} \label{EQ_107}
\widehat{\mathcal{P}_\phi}(z):=
\int_{-\infty}^{\infty} \mathfrak{S}_t^\sharp(z)\phi(t)  \, dt
~
\left(= 
\int_{-\infty}^{\infty} \overline{\mathfrak{S}_t(\bar{z})}\,\phi(t)  \, dt
\right)
\end{equation}
using \eqref{EQ_105}.  
Then $\widehat{\mathcal{P}_\phi}(z)$ belongs to $L^2(\R)$. 
\end{proposition}

Using the image of the composition $\widehat{\mathcal{P}_D}:=\widehat{\mathcal{P}} \circ D$ 
of the integral operator 
$\widehat{\mathcal{P}}$ and the differential operator 
\begin{equation} \label{EQ_108}
(D\psi)(t):=i\psi'(t), 
\end{equation}
we obtain the following equivalence condition for the RH.

\begin{theorem} \label{thm_1_2}
The RH is true if and only if the equality 
\begin{equation} \label{EQ_109}
\Vert \widehat{\mathcal{P}_{D\psi}} \Vert_{L^2(\R)}^2 = \pi \langle \psi, \psi \rangle_{W}
\end{equation}
holds for all $\psi \in C_c^\infty(\R)$. 
Furthermore, by choosing the test functions appropriately, 
if \eqref{EQ_109} holds for countably many choices of $\psi$'s, then the RH follows.
\end{theorem}
\begin{proof}
See Section \ref{section_4_3}. 
\end{proof}

Equation \eqref{EQ_109} is reformulated to the following simpler form. 

\begin{corollary} \label{cor_1_1} 
Define the subspace $V^\circ(0)$ of $L^2(\R)$ by 
\[
V^\circ(0) := \Bigl\{ \mathsf{F}^{-1}\widehat{\mathcal{P}_{D\psi}} \, \Bigl|\, 
\psi \in C_c^\infty(\R) \Bigr\}.  
\]
Then the RH is true if and only if the equality 
\begin{equation} \label{EQ_110}
2 \Vert \psi\Vert_{L^2(\R)}^2 = \langle \psi, \psi \rangle_{W} 
\end{equation}
holds for all  $\psi \in V^\circ(0)$. 
\end{corollary}
\begin{proof}
See Section \ref{section_4_3} and Theorem \ref{thm_5_2}. 
\end{proof}

The advantage of Theorem \ref{thm_1_2} and Corollary \ref{cor_1_1} 
is that it has turned the criterion of the RH 
from a set of inequalities like Weil's criterion into a set of equalities. 
It should also be noted that equations \eqref{EQ_109} and \eqref{EQ_110} can be expressed 
without zeros of $\xi(1/2-iz)$ by \eqref{EQ_105} and \eqref{EQ_106}. 
Furthermore, equations \eqref{EQ_109} and \eqref{EQ_110} claim that 
the nonnegativity of Weil's hermitian form 
is explained by the nonnegativity of the $L^2$-norm.
\medskip

In the following sections, first, 
in Section \ref{section_2}, 
we briefly review necessary notions such as 
screw functions, screw lines, the Hermite--Biehler class, de Branges spaces, and model spaces. 
Then, in Section \ref{section_3}, we state and prove unconditional results 
that we need to prove the main results.
Moreover, we unconditionally define two Hilbert spaces $\mathcal{H}_0$ and $\mathcal{K}_0$ that 
agree with the Hilbert spaces $\mathcal{H}_W$ and 
$\mathcal{K}(\Theta)$, respectively, under the RH. 

In Section \ref{section_4}, 
we show that $\mathfrak{S}_t(z)$ in \eqref{EQ_105} 
gives a screw line of the screw function corresponding to the Riemann zeta-function 
under the RH (Theorem \ref{thm_4_1}). 
Furthermore, we prove Theorem \ref{thm_1_2} and Corollary \ref{cor_1_1}. 
The strategy of the proof of Theorem \ref{thm_4_1} is basically the same as that of \cite[Theorem~1.1]{Su23}, with Proposition \ref{prop_401} playing an essential role in both cases.
To carry this out, the rewriting of \eqref{EQ_105} into \eqref{EQ_306}, 
prepared in Section~\ref{section_3} using Weil's explicit formula, 
corresponds to the transformation from (1.7) to (3.6) in \cite{Su23}, 
although the technical details of the calculations differ considerably.
On the other hand, the analytic or geometric meaning of the functions giving the norms 
was unclear in \cite{Su23}, 
whereas in the present article, these functions have a clear interpretation as a screw line. 
Furthermore, as an advantage of employing the screw line $\mathfrak{S}_t(z)$, 
we obtain Theorem \ref{thm_1_2}, for which no analog was obtained in \cite{Su23}.

In Section \ref{section_5}, we prove Theorem \ref{thm_1_1} in a more detailed form. 
In addition, we prove that $\mathcal{H}_0=\mathcal{H}_W$ 
and $\mathcal{K}_0=\mathcal{K}(\Theta_\xi)$ under the RH. 
Afterward, we explain that the Hilbert space $\mathcal{H}_W$ is one of the Hilbert--P{\'o}lya spaces 
in Section \ref{section_6}. 
Finally, we mention two special values of $\mathfrak{S}_t(z)$ in Section \ref{Section_4} 
as an appendix.

\section{Review on necessary notions} \label{section_2}

\subsection{Screw functions and screw lines}  \label{sec_2_1}
In this and the next part, we refer to \cite[Sections 5 and 12]{KrLa14}.  
See also its references for details. 
Following Kre\u{\i}n, 
we denote by $\mathcal{G}_\infty$ 
the space of all continuous functions $g(t)$ 
on $\R$ such that $g(-t) = \overline{g(t)}$ 
and the kernel 
\begin{equation} \label{EQ_201}
G_g(t,u):=g(t-u)-g(t)-g(-u)+g(0)
\end{equation}
is nonnegative definite on $\R$, that is,  
$\sum_{i,j=1}^{n} G_g(t_i,t_j) \,  \xi_i \overline{\xi_j} \,\geq\, 0$ 
for all $n \in \N$, $t_i \in \R$, and $\xi_i \in \C$ $(i = 1, 2, . . . , n)$. 
Functions belonging to $\mathcal{G}_\infty$ are called {\it screw functions} on $\R$. 

If an (even) real-valued function $g(t)$ is a screw function, 
then there exists a Hilbert space $\mathcal{H}$ 
and a continuous mapping $t \mapsto x(t)$ from $\R$ into $\mathcal{H}$ 
such that 
\[
\langle x(t+v)-x(v), x(u+v)-x(v) \rangle_{\mathcal{H}}
\]
is independent of $v \in \R$ for all $t,u \in \R$ 
and the equality 
$\langle x(t)-x(0), x(u)-x(0) \rangle_{\mathcal{H}} =G_{g}(t,u)$ holds.  
Therefore, $\Vert x(t)-x(0) \Vert_{\mathcal{H}}^2=-2g(t)$ 
if $g(0)=0$. 
A mapping $x:\R \to \mathcal{H}$ endowed with 
the translation-invariance described above is called a {\it screw line} 
for $g(t)$.

\subsection{Hilbert spaces associated with screw functions}

Each $g \in \mathcal{G}_\infty$ defines a nonnegative definite hermitian form on $\R$ by 
\begin{equation} \label{EQ_202}
\langle \phi_1,\phi_2 \rangle_{G_g}
:= \int_{-\infty}^{\infty}\int_{-\infty}^{\infty}G_g(t,u)\phi_1(u)\overline{\phi_2(t)} \, dudt. 
\end{equation}
According to \cite[Section 5]{KrLa14}, we denote by 
$\mathcal{L}(G_g)$ 
the space $C_0(\R)$ 
of all continuous and compactly supported functions $\phi$ on $\R$ 
such that $\widehat{\phi}(0)=0$ 
equipped with the hermitian inner product 
$\langle \cdot,\cdot \rangle_{G_g}$. 
We also denote by $\mathcal{H}(G_g)$ 
the completion of the factor space 
$\mathcal{L}(G_g)/\mathcal{L}^{\circ}(G_g)$, 
where 
$\mathcal{L}^{\circ}(G_g)=\{\phi \in \mathcal{L}(G_g)\,|\,\langle \phi,\phi \rangle_{G_g}=0\}$. 
Note that even if $\langle \cdot,\cdot \rangle_{G_g}$ 
is positive definite on $\mathcal{L}(G_g)$, 
that is, $\mathcal{L}^{\circ}(G_g)=\{0\}$, 
there possibly exists a sequence $(\phi_n)_n$ of  $\mathcal{L}(G_g)$ 
such that $\phi_n \to 0$ as $n \to \infty$  
with respect to $\langle \cdot,\cdot \rangle_{G_g}$. 
The completion $\mathcal{H}(G_g)$ is a space of equivalence classes of 
Cauchy sequences with respect to $\langle \cdot,\cdot \rangle_{G_g}$. 
Two Cauchy sequences are equivalent  
if their difference converges to zero with respect to $\langle \cdot,\cdot \rangle_{G_g}$. 
We denote by $[\phi] \in \mathcal{H}(G_g)$ 
the equivalence class represented by $\phi$.  
In general, elements of $\mathcal{H}(G_g)$ are not necessarily 
represented by functions unlike $\mathcal{H}_W$ 
(cf. \cite[Section 4.3]{KrLa14}).

Every $g \in \mathcal{G}_\infty$ admits a representation 
\begin{equation} \label{EQ_203} 
g(t) = g(0)+ibt + \int_{-\infty}^{\infty} \left( 
e^{i\lambda t}-1 - \frac{i\lambda t}{1+\lambda^2}
\right) \frac{d\tau(\lambda)}{\lambda^2}
\end{equation}
with $b \in \R$ and a measure $\tau$ on $\R$ 
such that $\int_{-\infty}^{\infty}d\tau(\lambda)/(1+\lambda^2)<\infty$ 
and vice versa. If $g(t)$ is real-valued, $b=0$. 
Without loss of generality, we suppose that $g(0)=0$. 

We define 
\[
\Phi_1(\phi,\lambda):=\int_{-\infty}^{\infty} 
\frac{e^{i\lambda x}-1}{\lambda} \,\phi(x) \, dx
= \frac{\widehat{\phi}(\lambda)-\widehat{\phi}(0)}{\lambda} 
= \frac{\widehat{\phi}(\lambda)}{\lambda} 
\]
for $\phi \in \mathcal{L}(G_g)$. 
Then, $\langle \phi_1, \phi_2 \rangle_{G_g}=\langle \Phi_1(\phi_1), \Phi_1(\phi_2) \rangle_{L^2(\tau)}$ 
for $\phi_1, \phi_2 \in \mathcal{L}(G_g)$ 
and $\Phi_1$ establishes an isomorphism between $\mathcal{H}(G_g)$ and $L^2(\tau)$.

\subsection{De Branges spaces} \label{sec_2_3}

In this part, we refer to \cite{SiTo15, Wo15}. 
See also those references for details. 
Let $H^2:=H^2(\C_+)=\mathsf{F}(L^2(0,\infty))$ 
be the Hardy space in the upper half-plane. 
As usual, 
we identify $H^2$ with a closed subspace of $L^2(\R)$ via boundary values. 
Then, the inner product of $H^2$ coincides with the standard inner product of $L^2(\R)$. 

The Hermite--Biehler class consists of entire functions $E$ satisfying 
$|E^\sharp(z)|<|E(z)|$ for all $z \in \C_+$.  
For each entire function $E$ belonging to the Hermite--Biehler class, 
the de Branges space $\mathcal{H}(E)$ is defined 
as a Hilbert space consisting of entire functions $F(z)$ 
such that both $F(z)/E(z)$ and $F^\sharp(z)/E(z)$ 
belong to $H^2$ and have the norm  
\begin{equation} \label{EQ_204}
\Vert F \Vert_{\mathcal{H}(E)} := \Vert F/E \Vert_{L^2(\R)}.
\end{equation}

The multiplication operator $\mathsf{M}$ by an independent variable 
is defined by 
$\mathfrak{D}(\mathsf{M})=\{ F(z) \in \mathcal{H}(E)\,|\, zF(z) \in \mathcal{H}(E)\}$ 
and $(\mathsf{M}F)(z)=zF(z)$ for $F \in \mathfrak{D}(\mathsf{M})$. 
The domain $\mathfrak{D}(\mathsf{M})$ is dense in $\mathcal{H}(E)$ 
if and only if 
\[
S_\theta(z) := \frac{i}{2}(e^{i\theta}E(z)-e^{-i\theta}E^\sharp(z))
\]
does not belong to $\mathcal{H}(E)$ for all $\theta \in [0,\pi)$ \cite[Theorem 11]{SiTo15}. 
The particular two $\theta$ cases are often written as   
$A(z):=-S_{\pi/2}(z)$ and  $B(z):=S_{0}(z)$. 

If $\mathfrak{D}(\mathsf{M})$ is dense in $\mathcal{H}(E)$, 
all self-adjoint extensions of $\mathsf{M}$ 
are parametrized by $\theta \in [0,\pi)$ and are described as follows. 
The domain of ${\mathsf M}_\theta$ is 
\begin{equation} \label{EQ_205}
{\mathfrak D}({\mathsf M}_\theta) 
= \left\{\left.
G(z) = \frac{S_\theta(w_0)F(z)-S_\theta(z)F(w_0)}{z-w_0} ~\right|~
F(z) \in {\mathcal H}(E)
\right\},
\end{equation}
and the operation is defined by 
\begin{equation} \label{EQ_206}
{\mathsf M}_\theta G(z) 
= z \, G(z) + F(w_0)S_\theta(z),
\end{equation}
where $w_0$ is a fixed complex number with $S_\theta(w_0)\not=0$ 
\cite[Propositions 4.6 and 6.1]{KW99}. 
The domain ${\mathfrak D}({\mathsf M}_\theta)$ is independent of the choice of the number $w_0$. 
For a fixed $\theta \in [0,\pi)$, we confirm that 
$G(z)=S_\theta(z)/(z-\gamma)$ belongs to $\mathfrak{D}(\mathsf{M}_\theta)$ 
by taking 
\[
F(z)=\frac{S_\theta(z)}{S_\theta(w_0)}\frac{\gamma-w_0}{z-\gamma}
\]
for every zero $\gamma$ of $S_\theta(z)$ 
and is an eigenfunction of $\mathsf{M}_\theta$ with the eigenvalue $\gamma$.  
Further, $\{S_\theta(z)/(z-\gamma)\,|\,S_\theta(\gamma)=0\}$ 
forms an orthogonal basis of $\mathcal{H}(E)$ \cite[Theorem 22]{dB68}.

\subsection{Model subspaces} \label{sec_2_4}

In this part, we refer to \cite[Section 2]{MaPo05}, \cite[Section 3.5]{Su20a} and \cite[Section 3.1]{Su23}. 
See also those references for details. 

Let $H^\infty=H^\infty(\C_+)$ be the space of all bounded analytic functions in $\C_+$. 
A function $\Theta \in H^\infty$  
is called an inner function in $\C_+$ 
if $\lim_{y \to 0+}|\Theta(x+iy)|=1$ for almost all $x \in \R$. 
For an inner function $\Theta$, 
a model space $\mathcal{K}(\Theta)$ is defined as the orthogonal complement 
$\mathcal{K}(\Theta)=H^2 \ominus \Theta H^2$ 
and has the alternative representation 
\begin{equation} \label{EQ_207}
\mathcal{K}(\Theta) = H^2 \cap \Theta \,\bar{H}^2, 
\end{equation} 
where $\Theta H^2 = \{ \Theta(z)F(z) \, |\, F \in H^2\}$ and $\bar{H}^2=H^2(\C_-)$ 
is the Hardy space in the lower half-plane.  
The model space $\mathcal{K}(\Theta)$ 
is a subspace of $L^2(\R)$ as a Hilbert space. 
In particular, the inner product of $\mathcal{K}(\Theta)$ matches 
that of $L^2(\R)$ on the real line. 

If an inner function $\Theta$ in $\C_+$ extends to a meromorphic function in $\C$, 
then it is called a meromorphic inner function in $\C_+$. 
For any meromorphic inner function $\Theta$, 
there exists $E$ of the Hermite--Biehler class 
such that $\Theta=E^\sharp/E$. 
The de Branges space $\mathcal{H}(E)$ 
is isometrically isomorphic to $\mathcal{K}(\Theta)$ by $F(z) \mapsto E(z)F(z)$. 
In particular, $\mathcal{H}(E) = E\,H^2 \cap E^\sharp \,\bar{H}^2$

For a meromorphic inner function $\Theta$, 
let $\mu_\Theta$ be the positive discrete measure on $\R$ 
supported on $\sigma(\Theta)=\{x\in \R\,|\,\Theta(x)=-1\}$ 
and 
\begin{equation} \label{EQ_208}
\mu_\Theta(x)=\frac{2\pi}{|\Theta'(x)|}.
\end{equation}
Then the restriction map $F \mapsto F|_{\sigma(\Theta)}$ 
is an isometric operator from $\mathcal{K}(\Theta)$ to $L^2(\mu_\Theta)$ 
\cite[Theorem 2.1]{MaPo05}. 
The isometric property of the map implies that the family of functions
\begin{equation} \label{EQ_209}
f_\gamma(z)=\sqrt{\frac{2}{\pi |\Theta'(\gamma)|}}\frac{1+\Theta(z)}{2(z-\gamma)}
=\sqrt{\frac{2}{\pi |\Theta'(\gamma)|}}\frac{A(z)}{(z-\gamma)E(z)}
\end{equation}
parametrized by all zeros $\gamma$ of $A(z)=-S_{\pi/2}(z)$ 
forms an orthonormal basis of $\mathcal{K}(\Theta)$ 
if $\mathfrak{D}(\mathsf{M})$ is dense in $\mathcal{H}(E)$.

\section{Unconditional results} \label{section_3}

Throughout this and later sections, 
we denote $E=E_\xi$ and $\Theta=\Theta_\xi=E_\xi^\sharp/E_\xi$ 
for functions defined in \eqref{EQ_103} and \eqref{EQ_104}, respectively. Otherwise, it is mentioned. 

\subsection{Expansion of $\mathfrak{P}_t(z)$ over the zeros} 
For the basic properties of the Riemann zeta-function, we refer to \cite{Tit86}. 
By the two functional equations $\xi(s)=\xi(1-s)$ and $\xi(s)=\xi^\sharp(s)$, 
if $\gamma$ belongs to the set of zeros $\Gamma$, 
then both $-\gamma$ and $\overline{\gamma}$ 
also belong to $\Gamma$ with the same multiplicity.  
On the other hand, $|\Im(\gamma)| < 1/2$ for every $\gamma \in \Gamma$, 
since all zeros of $\xi(s)$ lie in the strip $0 < \Re(s) < 1$. 
For $E(z)$ of \eqref{EQ_103}, we define 
\begin{equation}  \label{EQ_301}
A(z) := (E(z)+E^\sharp(z))/2 
\end{equation}
as in Section \ref{sec_2_3}. 
Then $A(z)=\xi(1/2-iz)$, because $E^\sharp(z)=\overline{E(\bar{z})}=\xi(1/2-iz)-\xi'(1/2-iz)$ 
by functional equations of $\xi(s)$. 
Therefore, the set $\Gamma$ coincides with the set of all zeros of both $A(z)$ and $1+\Theta(z)$. 
We define 
\begin{equation} \label{EQ_302} 
P_t(z) 
:= \sum_{\gamma \in \Gamma} m_\gamma \, \frac{e^{-i\gamma t}-1}{\gamma}
\cdot
\frac{1}{z-\gamma}
\end{equation}
for nonnegative $t$.  
For negative $t$, we set $P_{t}(z):=P_{-t}(z)$. 
The series on the right-hand side of \eqref{EQ_302} converges absolutely and uniformly 
on every compact subset of $\C\setminus\Gamma$, 
since $\sum_{\gamma \in \Gamma}m_\gamma|\gamma|^{-1-\delta}<\infty$ for any $\delta>0$, 
because $A(z)$ is an entire function of order one. 
Therefore, $P_t(z)$ is a meromorphic function on $\C$ 
with $\Gamma$ as the set of all poles.

\begin{proposition} \label{prop_301} 
Let $\mathfrak{P}_t(z)$ and $P_t(z)$ be meromorphic functions 
defined by \eqref{EQ_106} and \eqref{EQ_302}, respectively. 
Then, both coincide.
\end{proposition}
\begin{proof} For $t \geq 0$ and $z \in \C_+$, we define 
\[
\phi_{z,t}(x) = 
\begin{cases}
~(iz)^{-1} \,e^{izx} (e^{-izt}-1), &  t < x, \\
~(iz)^{-1} \,e^{izx} (e^{-izx}-1), & 0 \leq x \leq t, \\
~0, & x<0.   
\end{cases}
\]
The main tool for the proof is the Weil explicit formula
\begin{equation} \label{EQ_303}
\aligned 
\lim_{X \to \infty} & \sum_{{\gamma \in \Gamma}\atop{|\gamma|\leq X}} m_\gamma
\int_{-\infty}^{\infty} \phi(x) \, e^{-i\gamma x} \, dx \\
& = \int_{-\infty}^{\infty} \phi(x) (e^{x/2}+e^{-x/2}) dx 
 - \sum_{n=1}^{\infty} \frac{\Lambda(n)}{\sqrt{n}} \phi(\log n) 
- \sum_{n=1}^{\infty} \frac{\Lambda(n)}{\sqrt{n}} \phi(-\log n)  \\
& \quad 
- (\log 4\pi + \gamma_0) \phi(0)
- \int_{0}^{\infty} \left\{ \phi(x) + \phi(-x)-2e^{-x/2} \phi(0)\right\} \frac{e^{x/2}dx}{e^{x}-e^{-x}}
\endaligned 
\end{equation}
which is obtained from the explicit formula in \cite[p. 186]{Bo01} 
by taking $\phi(x) = e^{x/2}f(e^x)$ for test functions $f(t)$ in that formula
with the conditions for $f(t)$ in \cite[Section 3]{BoLa99}, 
where $\gamma_0$ is the Euler--Mascheroni constant. 
(Note that the formula in \cite{BoLa99} has two typographical errors 
in the second line of the right-hand side.) 

As is easily seen, Weil's explicit formula can be applied to $\phi(x)=\phi_{z,t}(x)$. 
We have  
\[
\int_{-\infty}^{\infty} \phi_{z,t}(x) \, e^{-i\gamma x} \, dx 
= \frac{e^{-i\gamma t}-1}{\gamma}
\cdot
\frac{1}{z-\gamma} 
\quad \text{when $\Im(z)>\Im(\gamma)$}. 
\]
Therefore, the left-hand side of Weil's explicit formula for $\phi_{z,t}(x)$ 
gives $P_t(z)$ of \eqref{EQ_302} when $\Im(z)>1/2$. 
Hence, if it is shown that the right-hand side 
is equal to $\mathfrak{P}_t(z)$ for $\Im(z)>1/2$, 
then the conclusion of the proposition follows by analytic continuation. 

It is easy to verify
\[
\int_{-\infty}^{\infty} \phi_{z,t}(x) (e^{x/2}+e^{-x/2}) dx 
= \frac{4(e^{t/2}-1)}{1+2iz} + \frac{4(e^{-t/2}-1)}{1-2iz}  
\]
and 
\[
\aligned 
\sum_{n=1}^{\infty} \frac{\Lambda(n)}{\sqrt{n}}  \phi_{z,t}(\log n) 
& = \frac{1}{iz} \sum_{ n \leq e^t} \frac{\Lambda(n)}{\sqrt{n}} (1-n^{iz})
 + \frac{e^{-izt}-1}{iz}\sum_{t < \log n} \frac{\Lambda(n)}{n^{1/2-iz}} 
\\
& = - \sum_{n \leq e^t} \frac{\Lambda(n)}{\sqrt{n}} \frac{e^{-iz(t-\log n)}-1}{iz}
 - \frac{e^{-izt}-1}{iz}\frac{\zeta'}{\zeta}\left( \frac{1}{2}-iz \right), 
\\
\sum_{n=1}^{\infty} \frac{\Lambda(n)}{\sqrt{n}}  \phi_{z,t}(-\log n) 
& = 0, \qquad \phi_{z,t}(0) =0 
\endaligned 
\]
for $\Im(z)>1/2$ by direct calculation. 

Therefore, the remaining task is to calculate the fifth term on the right-hand side. 
We split it into $\int_{t}^{\infty}$ and $\int_{0}^{t}$. 
For the first integral, 
\[
\aligned 
\int_{t}^{\infty} & \left\{ \phi_{z,t}(x) + \phi_{z,t}(-x)-2e^{-x/2} \phi_{z,t}(0)\right\} 
\frac{e^{x/2}dx}{e^{x}-e^{-x}} \\
& = \frac{e^{-izt}-1}{iz}  
\int_{t}^{\infty} e^{izx} \frac{e^{x/2}dx}{e^{x}-e^{-x}} 
 = \frac{e^{-izt}-1}{iz}  
\int_{t}^{\infty} e^{izx} \,e^{-x/2} \sum_{n=0}^{\infty}e^{-2nx} \, dx \\
& = \frac{e^{-izt}-1}{2iz} 
e^{-t(\frac{1}{2}-iz)}
\sum_{n=0}^{\infty} \frac{e^{-2nt}}{n+\frac{1}{2}(\frac{1}{2}-iz)} 
 = \frac{e^{-izt}-1}{2iz} 
e^{-t(\frac{1}{2}-iz)}
\Phi(e^{-2t},1,\tfrac{1}{2}(\tfrac{1}{2}-iz)). 
\endaligned 
\]
For the second integral, 
\[
\aligned 
\int_{0}^{t} & \left\{ \phi_{z,t}(x) + \phi_{z,t}(-x)-2e^{-x/2} \phi_{z,t}(0)\right\} \frac{e^{x/2}dx}{e^{x}-e^{-x}} \\
& = -\frac{1}{iz} \int_{0}^{t} (e^{izx}-1) \frac{e^{x/2}dx}{e^{x}-e^{-x}} 
= -\frac{1}{iz} \int_{0}^{t} (e^{izx}-1) \,e^{-x/2} \sum_{n=0}^{\infty}e^{-2nx} \, dx. 
\endaligned 
\]
To handle the right-hand side, we calculate as 
\[
\aligned
\int_{0}^{t} (e^{izx}-1) & \,e^{-x/2} \sum_{n=0}^N e^{-2nx} \, dx \\
& = \frac{1}{2} \sum_{n=0}^N
\left[
\frac{1-e^{-2t(n+\frac{1}{2}(\frac{1}{2}-iz))}}{n+\frac{1}{2}(\frac{1}{2}-iz)}
-
\frac{1-e^{-2t(n+\frac{1}{4})}}{n+\frac{1}{4}}
\right] \\
& = -\frac{1}{2} e^{-t(\frac{1}{2}-iz)}\sum_{n=0}^N
\frac{e^{-2tn}}{n+\frac{1}{2}(\frac{1}{2}-iz)}
+ \frac{1}{2} e^{-t/2} \sum_{n=0}^N
\frac{e^{-2tn}}{n+\frac{1}{4}}
 \\
& \quad +\frac{1}{2} \sum_{n=0}^N
\left[
\frac{1}{n+\frac{1}{2}(\frac{1}{2}-iz)}
-
\frac{1}{n+\frac{1}{4}}
\right] \\
& = -\frac{1}{2} e^{-t(\frac{1}{2}-iz)}
\Phi(e^{-2t},1,\tfrac{1}{2}(\tfrac{1}{2}-iz))
+\frac{1}{2} e^{-t/2} \Phi(e^{-2t},1,\tfrac{1}{4})
 \\
& \quad - \frac{1}{2} 
\left[
\frac{\Gamma'}{\Gamma}\left(\frac{1}{4}-\frac{iz}{2}\right)
-
\frac{\Gamma'}{\Gamma}\left(\frac{1}{4}\right)
\right] +O(e^{-2Nt}) + O(N^{-1}) 
\endaligned 
\]
using the well-known series expansion
\begin{equation} \label{EQ_304} 
\frac{\Gamma'}{\Gamma}(w) = -\gamma_0 - \sum_{n=0}^{\infty}
\left( \frac{1}{w+n} - \frac{1}{n+1} \right), 
\end{equation}
where the implied constant depends on $t$ and $z$. 
Therefore, we obtain
\[
\aligned 
\int_{0}^{t} & \left\{ \phi_{z,t}(x) + \phi_{z,t}(-x)-2e^{-x/2} \phi_{z,t}(0)\right\} \frac{e^{x/2}dx}{e^{x}-e^{-x}} \\
& = \frac{1}{2iz} e^{-t(\frac{1}{2}-iz)}
\Phi(e^{-2t},1,\tfrac{1}{2}(\tfrac{1}{2}-iz))
-\frac{1}{2iz} e^{-t/2} \Phi(e^{-2t},1,\tfrac{1}{4})
 \\
& \quad +\frac{1}{2iz} 
\left[
\frac{\Gamma'}{\Gamma}\left(\frac{1}{4}-\frac{iz}{2}\right)
-
\frac{\Gamma'}{\Gamma}\left(\frac{1}{4}\right)
\right].
\endaligned 
\]
Combining the results for $\int_{t}^{\infty}$ and $\int_{0}^{t}$, 
\[
\aligned 
\int_{0}^{\infty} & 
\left\{ \phi_{z,t}(x) + \phi_{z,t}(-x)-2e^{-x/2} \phi_{z,t}(0)\right\} \frac{e^{x/2}dx}{e^{x}-e^{-x}} \\
& = \frac{1}{2iz} 
e^{-t/2} \Bigl[
\Phi(e^{-2t},1,\tfrac{1}{2}(\tfrac{1}{2}-iz)) 
-\Phi(e^{-2t},1,\tfrac{1}{4}) 
\Bigr]
 \\
& \quad +\frac{1}{2iz} 
\left[
\frac{\Gamma'}{\Gamma}\left(\frac{1}{4}-\frac{iz}{2}\right)
-
\frac{\Gamma'}{\Gamma}\left(\frac{1}{4}\right)
\right]. 
\endaligned 
\]
From the calculation of the five terms on the right-hand side above,
we conclude that the right-hand side of the Weil explicit formula 
for $\phi_{z,t}(x)$ equals \eqref{EQ_106}.
\end{proof}

\subsection{Proof of Proposition \ref{prop_101}} \label{sec_3_2}

%
We have $|\Theta(z)|=1$ for every $z \in \R$ by definition. 
In fact, zeros of $E(z)$ in the denominator cancel out in the numerator $E^\sharp(z)$, 
even if they exist. 
Further, $\mathfrak{P}_t(z)$ has poles of order one at $\gamma \in \Gamma$, 
but $\mathfrak{S}_t(z)$ is holomorphic there, since 
$(1+\Theta(z))/2 = A(z)/E(z) = A(z)/(A(z)+iA'(z)) = (z-\gamma)(-i/m_\gamma+o(1)) 
$ near $z=\gamma$ by direct calculation.  
Hence, $\mathfrak{S}_t(z)$ is bounded and holomorphic on the real line 
by \eqref{EQ_105}, \eqref{EQ_302}, and Proposition \ref{prop_301}. 
On the other hand, in the horizontal strip $|\Im(z)|\leq 1/2$, 
we have the well-known estimate $(\Gamma'/\Gamma)(1/4+iz/2) \ll \log |z|$ 
and 
\[
\frac{\zeta'}{\zeta}\left(\frac{1}{2}-iz \right)
= \sum_{|\Re(z)-\gamma| \leq 1} \frac{i}{z-\gamma}+O(\log|z|)
\]
by \cite[Theorem 9.6 (A)]{Tit86}. 
In both estimates, implied constants are uniform in $|\Im(z)|\leq 1/2$. 
The number of zeros $\gamma \in \Gamma$ satisfying $|\Re(z)-\gamma| \leq 1$ 
is $O(\log |z|)$ counting with multiplicity by \cite[Theorem 9.2]{Tit86}. 
Therefore, 
$\mathfrak{S}_t(z) \ll |z|^{-1}\log |z|$ as $|z| \to \infty$ 
with an implied constant depending on a compact set of $t$ by \eqref{EQ_106}. 
Hence $\mathfrak{S}_t(z)$ belongs to $L^2(\R)$ 
and the norm is uniformly bounded on a compact set of $t$. 
\hfill $\Box$ 

\subsection{Two special Hilbert spaces} \label{section_3_3}

We first introduce the set of meromorphic functions 
\begin{equation} \label{EQ_305}
F_\gamma(z) :=
\sqrt{\frac{m_\gamma}{\pi}} \frac{i(1+\Theta(z))}{2(z-\gamma)}, \quad \gamma \in \Gamma. 
\end{equation}
Then, we have 
\begin{equation} \label{EQ_306} 
\mathfrak{S}_t(z)
= 
\sum_{\gamma \in \Gamma} \sqrt{\pi m_\gamma} \, \frac{e^{-i\gamma t}-1}{\gamma} \, F_\gamma^\sharp(z) 
\end{equation}
by Proposition \ref{prop_301}. 
Therefore, 
\begin{equation} \label{EQ_307} 
\widehat{\mathcal{P}_\phi}(z)
= 
\sum_{\gamma \in \Gamma} \sqrt{\pi m_\gamma} \, 
\frac{\widehat{\phi}(\gamma)-\widehat{\phi}(0)}{\gamma} \, F_\gamma(z) 
\end{equation}
for any $\phi \in C_c^\infty(\R)$ by definition \eqref{EQ_107} 
and the symmetry $\gamma \mapsto \bar{\gamma}$ of $\Gamma$ 
with $m_\gamma=m_{\bar{\gamma}}$. 
This implies 
\begin{equation} \label{EQ_308} 
\widehat{\mathcal{P}_{D\psi}}(z)
= 
\sum_{\gamma \in \Gamma} \sqrt{\pi m_\gamma} \, 
\widehat{\psi}(\gamma) \, F_\gamma(z) 
\end{equation}
for any $\psi \in C_c^\infty(\R)$, 
since 
$(\widehat{D\psi}(z)-\widehat{D\psi}(0))/z
=\widehat{D\psi}(z)/z=\widehat{\psi}(z)$ 
for $D$ in \eqref{EQ_108}. 
\medskip

On the other hand, 
we define the norm $\Vert~\Vert_0$ on $C_c^\infty(\R)$ by 
\begin{equation} \label{EQ_309}
\Vert \psi \Vert_0 := \frac{1}{\sqrt{\pi}} \, \Vert  \widehat{\mathcal{P}_{D\psi}} \Vert_{L^2(\R)}, 
\quad \psi \in C_c^\infty(\R)
\end{equation} 
based on Proposition \ref{prop_102}. Then, we have the following.

\begin{lemma} \label{lem_301}
Equation \eqref{EQ_309} defines a norm on $C_c^\infty(\R)$. 
\end{lemma}
\begin{proof}
We obtain $\Vert \psi_1+\psi_2 \Vert_0 \leq \Vert \psi_1 \Vert_0 + \Vert \psi_2 \Vert_0$ 
and $\Vert k \psi \Vert_0 = |k| \Vert \psi \Vert_0$ 
for $\psi_1, \psi_2, \psi \in C_c^\infty(\R)$ and $k \in \C$ 
by the obvious linearity  of $\widehat{\mathcal{P}_D}$. 
Therefore, 
the proof is completed if it is shown that $\Vert \psi \Vert_0=0$ implies $\psi=0$. 
If $\Vert \psi \Vert_0=0$,  
the image $\widehat{\mathcal{P}_{D\psi}}(z)$ is identically zero. 
The latter means that $\widehat{\psi}(\gamma)=0$ for all $\gamma \in \Gamma$, 
because, if not, there must exist a sequence $(c_\gamma)_{\gamma \in \Gamma}$ 
such that $\sum_{\gamma \in \Gamma} c_\gamma(z-\gamma)^{-1}$ is identically zero on $\C$ 
by \eqref{EQ_305} and \eqref{EQ_308}, but it is impossible. 
If $\widehat{\psi}(\gamma)=0$ for all $\gamma \in \Gamma$, 
it implies that $\psi$ is identically zero by \cite[Lemma 2.1]{Su22}. 
\end{proof}

By Lemma \ref{lem_301}, 
we can complete the space $C_c^\infty(\R)$ with respect to $\Vert~\Vert_0$. 
We denote the completion by $\mathcal{H}_0$. 
On the other hand, 
we denote the $L^2$-closure of 
the image $\widehat{\mathcal{P}_{D}}(C_c^\infty(\R))$ in $L^2(\R)$ by $\mathcal{K}_0$. 
Then, two Hilbert spaces $\mathcal{H}_0$ and $\mathcal{K}_0$ 
are isometrically isomorphic up to a constant multiple. 
The map $\widehat{\mathcal{P}_D}$ from $C_c^\infty(\R)$ to 
$\widehat{\mathcal{P}_{D}}(C_c^\infty(\R)) \subset L^2(\R)$ 
extends to the map from $\mathcal{H}_0$ to $\mathcal{K}_0$ by \eqref{EQ_309}. 
As proved in Theorem \ref{thm_5_2} below, 
$\mathcal{H}_0=\mathcal{H}_W$ and $\mathcal{K}_0=\mathcal{K}(\Theta)$ 
under the RH. 

\section{A screw line of the Riemann zeta-function} \label{section_4}

\subsection{A special orthonormal basis} 
Assuming the RH is true, 
$E=E_\xi$ belongs to the Hermite--Biehler class \cite[Theorem 1]{La06}, 
and thus $\Theta=\Theta_\xi$ is a meromorphic inner function. 
Therefore, they define the de Branges space $\mathcal{H}(E)$ 
and the model space $\mathcal{K}(\Theta)$, respectively. 
We need the following result for the later discussion.

\begin{proposition} \label{prop_401}
Assume that the RH is true. Then, the family \eqref{EQ_305} 
forms an orthonormal basis of the Hilbert space $\mathcal{K}(\Theta)$. 
Furthermore, 
\begin{equation} \label{EQ_401}
\frac{\Theta'(\gamma)}{2} 
= -\frac{i}{m_\gamma}
\end{equation}
and 
\begin{equation} \label{EQ_402} 
F_\gamma(\gamma) 
=  \frac{1}{\sqrt{m_\gamma \pi}}, 
\qquad F_\gamma(\gamma')=0 \quad \text{for every} \quad
\gamma \in \Gamma, ~
 \gamma' \in \Gamma\setminus\{\gamma\}.
\end{equation}

\end{proposition}
\begin{proof} 
See \cite[Proposition 3.2]{Su23} and its proof. 
\end{proof}

\subsection{Screw line of the Riemann zeta-function} \label{section_4_2}

We define the even real-valued function $g_\xi(t)$ on the real line by 
\begin{equation} \label{EQ_403}
\aligned 
g_\xi(t)
& := -4(e^{t/2}+e^{-t/2}-2)  + \sum_{n \leq e^t} \frac{\Lambda(n)}{\sqrt{n}}(t-\log n)
\\
&\quad  
- 
\frac{t}{2}\left[ \frac{\Gamma'}{\Gamma}\left(\frac{1}{4}\right) - \log \pi \right] 
- 
\frac{1}{4}\left( \Phi(1,2,1/4) - e^{-t/2}\Phi(e^{-2t},2,1/4) \right)
\endaligned 
\end{equation}
for nonnegative $t$. We easily obtain $g_\xi(0)=0$. 
Then, $g_\xi(t)$ is a screw function on $\R$ 
under the RH as stated in \cite[Theorem 1.2]{Su22}. 
One of the screw lines corresponding to $g_\xi(t)$ can be constructed as follows. 

Let  $\tau_\xi$ be the nonnegative measure 
representing $g_\xi(t)$ as in \eqref{EQ_203} under the RH. 
Then the Hilbert space $\mathcal{H}=L^2(\tau_\xi)$ and 
the mapping $t \mapsto x(t):=(e^{it\gamma}-1)/\gamma$ 
provide a screw line satisfying $\Vert x(t)-x(0) \Vert_{\mathcal{H}}^2=-2g_\xi(t)$ 
\cite[Section 12]{KrLa14}. 
This spectral construction of the screw line is important and useful in analysis, 
but it is of limited use for studying the nontrivial zeros of $\zeta(s)$ without assuming the RH. 
In the following, we show that $\mathfrak{S}_t$ gives a screw line of $g_\xi(t)$. 
In contrast to the spectral screw line above, this screw line can be used 
to study $\mathcal{H}_W$, as will be done later.

\begin{theorem} \label{thm_4_1} 
Assume the RH is true and let $g(t)=g_\xi(t)$.  
Then, the mapping 
$t \mapsto \pi^{-1/2} \mathfrak{S_t}(z)$ from $\R$ to $L^2(\R)$ 
is a screw line of $g(t)$. 
That is, 
\begin{equation} \label{EQ_404}
\frac{1}{\pi} \langle \mathfrak{S}_{t}, \mathfrak{S}_{u} \rangle_{L^2(\R)}
= G_g(t,u)
\end{equation}
holds for $t, u \in \R$.  
\end{theorem}
\begin{proof}
The sum of coefficients on the right-hand side of \eqref{EQ_306} 
is convergent in $L^2$-sense:
\[
\sum_{\gamma \in \Gamma}  \left|\sqrt{\pi m_\gamma} 
\frac{e^{-i\gamma t}-1}{\gamma}
 \right|^2 
\leq \pi \sum_{\gamma \in \Gamma} 
\frac{m_\gamma}{|\gamma|^2} < \infty. 
\]
Therefore,  applying Proposition \ref{prop_401} 
to $\mathfrak{S}_t(z)$ via formula \eqref{EQ_306}, 
we find that it belongs to the subspace $\mathcal{K}(\Theta)$ of $L^2(\R)$ and 
\begin{equation} \label{EQ_405}
\frac{1}{\pi} \langle \mathfrak{S}_{t+v}-\mathfrak{S}_{v}, \mathfrak{S}_{u+v}-\mathfrak{S}_{v} \rangle_{L^2(\R)}
= 
\sum_{\gamma \in \Gamma} m_\gamma \, 
\frac{e^{-i\gamma t}-1}{\gamma}\cdot\frac{e^{i\gamma u}-1}{\gamma}
\end{equation}
holds. The right-hand side is equal to $G_g(t,u)$ by the formula 
\begin{equation} \label{EQ_409}
G_g(t,u) 
 = \sum_{\gamma \in \Gamma} m_\gamma \frac{(e^{i\gamma t}-1)(e^{-i\gamma u}-1)}{\gamma^2}
\end{equation}
given in \cite[(1.9)]{Su22} and the symmetry $\gamma \mapsto -\gamma$ of $\Gamma$ 
with $m_\gamma = m_{-\gamma}$.  
Hence, $\pi^{-1/2}\mathfrak{S}_t:\R \to L^2(\R)$ 
is a screw line of $g(t)$ under the RH. 

We find that $\mathfrak{S}_0(z)$ is identically zero by \eqref{EQ_105} and \eqref{EQ_106}, 
since 
\[
\lim_{t  \to 0} \left(\Phi(e^{-2t},1,\tfrac{1}{4}) -\Phi(e^{-2t},1,\tfrac{1}{2}(\tfrac{1}{2}-iz))
\right) = 
-
\frac{\Gamma'}{\Gamma}\left(\frac{1}{4}\right)
+
\frac{\Gamma'}{\Gamma}\left(\frac{1}{2}\left(\frac{1}{2}-iz\right)\right) 
\]
by \eqref{EQ_208}. Therefore, by taking $v=0$ in \eqref{EQ_405}, we obtain \eqref{EQ_404}. 
\end{proof}

The following immediately follows from Theorem \ref{thm_4_1}. 

\begin{corollary} \label{cor_401}
The RH is true if and only if the equality 
\begin{equation} \label{EQ_406}
\frac{1}{2\pi} \Vert \mathfrak{S}_t \Vert_{L^2(\R)}^2 = - g(t)
\end{equation}
holds for all $t \geq t_0$ for some $t_0  \geq 0$. 
\end{corollary}
\begin{proof}
Assuming the RH, we obtain \eqref{EQ_406} by taking $u=t$ in \eqref{EQ_404}, 
since $G_g(t,t) = -2g(t)$ by \eqref{EQ_201} and $g(0)=0$. 
Conversely, 
we suppose that equality \eqref{EQ_406} holds for all $t \geq t_0$.  
Then $-g(t)$ is nonnegative on $[t_0,\infty)$, which implies 
that  the RH is true by \cite[Theorems 1.7 and 11.1]{Su22}. 
\end{proof}

\subsection{Proof of Theorem \ref{thm_1_2}} \label{section_4_3}

Theorem \ref{thm_1_2} is a corollary of the following result. 

\begin{theorem} \label{thm_4_2} 
Let $g(t)=g_\xi(t)$. 
The RH is true if and only if the equality 
\begin{equation} \label{EQ_407}
\Vert \widehat{\mathcal{P}_\phi} \Vert_{L^2(\R)}^2 = \pi \langle \phi, \phi \rangle_{G_g}
\end{equation}
holds for all $\phi \in C_c^\infty(\R)$ satisfying $\widehat{\phi}(0)=0$. 
If the RH is true, equality \eqref{EQ_407} holds for all $\phi \in C_c^\infty(\R)$. 
\end{theorem}
\begin{proof}
First, we prove \eqref{EQ_407} assuming the RH holds.   
We have
\begin{equation} \label{EQ_408} 
\Vert \widehat{\mathcal{P}_\phi} \Vert_{L^2(\R)}^2 
= \pi \sum_{\gamma \in \Gamma} m_\gamma
\left\vert \frac{\widehat{\phi}(\gamma)-\widehat{\phi}(0)}{\gamma} \right\vert^2
\end{equation}
by \eqref{EQ_307} and Proposition \ref{prop_401}. 
Applying \eqref{EQ_409} to \eqref{EQ_202} 
and noting the symmetry $\gamma \mapsto -\gamma$ 
of $\Gamma$ with $m_\gamma=m_{-\gamma}$, 
we find that the right-hand side of \eqref{EQ_408} equals 
$\pi \langle \phi, \phi \rangle_{G_g}$. 

Conversely, we prove that the RH is true assuming equality \eqref{EQ_407}. 
We show that a contradiction arises if the RH is false.  
We take a nonreal $\gamma_0 \in \Gamma$. 
For any $\epsilon>0$, there exists $\psi_1$, $\psi_2\in C_c^\infty(\R)$ 
such that $\widehat{\psi_1}(-\gamma_0)=i$, 
$\widehat{\psi_2}(-\overline{\gamma_0})=-i$, 
$|\widehat{\psi_1}(-\gamma)| \leq \epsilon |\gamma_0-\gamma|^{-1-\delta}$ 
for every $\gamma \in \Gamma\setminus\{\gamma_0\}$, 
and  
$|\widehat{\psi_2}(-\gamma)| \leq \epsilon |\overline{\gamma_0}-\gamma|^{-1-\delta}$ 
for every $\gamma \in \Gamma\setminus\{\overline{\gamma_0}\}$ by \cite[Lemma 1]{Yo92}. 
We define $\psi:=\psi_1+\psi_2 \,(\not=0)$ and set $\phi:=D\psi$.  
Then, $\widehat{\phi}(0)=0$ by definition, and 
$
\langle \phi, \phi \rangle_{G_g} 
=
\langle \psi, \psi \rangle_W 
$ 
holds by the relation 
\begin{equation} \label{EQ_410}
\langle D\psi_1, D\psi_2 \rangle_{G_g} = \langle \psi_1, \psi_2 \rangle_W
\end{equation}
%
in \cite[Proposition 3.1]{Su22}. The right-hand side equals 
$
\sum_{\gamma \in \Gamma} m_\gamma 
\widehat{\psi}(-\gamma)(\widehat{\psi})^\sharp(-\gamma) 
= -m_{\gamma_0} + O(\epsilon)
$,  
since $\sum_{\gamma \in \Gamma} m_\gamma |\gamma|^{-1-\delta}<\infty$. 
Therefore, $\langle \phi, \phi \rangle_{G_g}$ is negative for a sufficiently small $\epsilon>0$, 
but it contradicts the nonnegativity that follows from \eqref{EQ_407}. 
\end{proof}

\begin{proof}[Proof of Theorem \ref{thm_1_2}]

The conclusion follows 
from Theorem \ref{thm_4_2} and  the relation \eqref{EQ_410} 
of hermitian forms, 
since the differential operator $D$ in \eqref{EQ_108} gives a bijection from $C_c^\infty(\R)$ 
to the subspace 
$C_0^\infty(\R) \subset C_c^\infty(\R)$ consisting of functions $\phi$ 
with $\widehat{\phi}(0)=0$. 
\end{proof}

\begin{proof}[Proof of Corollary \ref{cor_1_1}]
The RH is true if \eqref{EQ_110} holds 
by the same argument as the second half of the proof of Theorem \ref{thm_4_2}. 
Therefore, we prove \eqref{EQ_110} assuming the RH. 

Let $\psi \in V^\circ(0)$. Then $\widehat{\psi}(z)=\widehat{\mathcal{P}_{D\psi_0}}(z)$ 
for some $\psi_0 \in C_c^\infty(\R)$ by definition. 
Therefore, 
$ 
\widehat{\psi}(z)
= 
\sum_{\gamma \in \Gamma} \sqrt{\pi m_\gamma} \, 
\widehat{\psi_0}(\gamma) F_\gamma(z) 
$ 
by \eqref{EQ_308}. 
The equality shows that $\widehat{\psi}(z)$ is a continuous function of $z \in \R$ 
by the uniform convergence of the right-hand side on a compact set of $z$.  
Taking $z=\gamma$ in this equality, 
we have 
$\widehat{\psi}(\gamma) = \widehat{\psi_0}(\gamma)$ 
by \eqref{EQ_402}. 
Therefore, $\langle \psi, \psi \rangle_W$ is defined 
and satisfies $\langle \psi, \psi \rangle_W = \langle \psi_0, \psi_0 \rangle_W$. 
The right-hand side is equal to 
$\Vert \widehat{\psi_0} \Vert_{L^2(\R)}^2 = 2\pi \Vert \psi_0 \Vert_{L^2(\R)}^2$
 by \eqref{EQ_109} and Plancherel's identity. 
The same argument works if we start with $\psi_0 \in C_c^\infty(\R)$. 
Hence, we obtain \eqref{EQ_110}. 
\end{proof}

Using \eqref{EQ_309}, Theorem \ref{thm_1_2} is stated as follows.

\begin{theorem} \label{thm_4_3}
The RH is true if and only if the equality 
\begin{equation} \label{EQ_411}
\Vert \psi\Vert_0^2 = \langle \psi, \psi \rangle_{W}
\end{equation}
holds for all $\psi \in C_c^\infty(\R)$. 
\end{theorem}

\noindent
Equality \eqref{EQ_411} leads to Theorem \ref{thm_5_2} below.
\medskip

For $n \in \Z_{>0}$, we define 
\[
g_n(x) := e^{-x/2}\sum_{j=1}^{n}\binom{n}{j} \frac{(-x)^{j-1}}{(j-1)!} \quad (x >0), \quad
g_n(0) := \frac{n}{2}, \quad 
g_n(x) :=0 \quad (x<0). 
\]
Then, the RH holds if $\langle g_n, g_n \rangle_{W}\geq 0$ for all $n \in \Z_{>0}$ 
by Bombieri and Lagarias \cite[Section 4]{BoLa99}. 
Therefore, we obtain the following. 

\begin{corollary} \label{cor_4_2}
The RH holds if \eqref{EQ_411} holds for all $g_n$ ($n \in \Z_{>0}$). 
\end{corollary}

\section{Proof of Theorem \ref{thm_1_1} and its refinement} \label{section_5}

Throughout this section, 
we assume that the RH is true 
and denote $E=E_\xi$, $\Theta=\Theta_\xi=E_\xi^\sharp/E_\xi$ as before, 
and denote $g=g_\xi$. 
Therefore, 
$E$ belongs to the Hermite--Biehler class, 
$\Theta$ is a meromorphic inner function in $\C_+$, 
and $g$ belongs to the class of screw functions $\mathcal{G}_\infty$. 
\medskip

For use in the proof of Theorem \ref{thm_1_1}  and its refinement, 
we introduce  the operator $\mathsf{K}$ acting on $L^2(\R)$ 
by 
\begin{equation} \label{EQ_501}
\mathsf{K}:=\mathsf{F}^{-1}\mathsf{M}_{\Theta}\mathsf{J}\mathsf{F} 
\end{equation}
with  
\[
(\mathsf{M}_{\Theta}F)(z):=\Theta(z)F(z) 
\quad \text{and} \quad (\mathsf{J}F)(z):=F^\sharp(z). 
\]
The Fourier transform $\mathsf{F}$, 
the multiplication operator $\mathsf{M}_{\Theta}$, 
and the involution $\mathsf{J}$ are 
defined for functions of a complex variable, 
and the latter two are isometries on $L^2(\R)$. 
The Fourier transform $\mathsf{F}$ is an isometry up to a constant factor.  
Therefore, $\mathsf{K}$ is isometric on $L^2(\R)$. 
Further, $\mathsf{K}$ is invertible by $\mathsf{K}^2={\rm id}$. 
By definition, $\mathsf{K}$ is {\it not} $\C$-linear but  $\R$-linear and conjugate linear. 
Using the isometric involution $\mathsf{K}$, we define
\begin{equation} \label{EQ_502}
V(t) := L^2(t,\infty) \cap \mathsf{K}L^2(t,\infty)
\end{equation}
and 
\[
\mathcal{H}_W(t):=\{\,[\psi]~|~\psi \in V(t)\,\}
\]
for $t \geq 0$. 
The set of subspaces $V(t)$ of $L^2(\R)$ are clearly totally ordered by the set-theoretical inclusion. 

First, Theorem \ref{thm_1_1} is shown using $V(t)$ for $t=0$, 
and it is refined using general $t \geq 0$.

\begin{lemma} \label{lem_5_1} 
Let $V(0) = L^2(0,\infty) \cap \mathsf{K}L^2(0,\infty)$. 
Then, we have 
$\mathcal{K}(\Theta)=\mathsf{F}(V(0))$, 
and hence $\mathcal{H}(E)=E \mathsf{F}(V(0))$
$=\{E(z)\widehat{\psi}(z)\,|\, \psi \in V(0)\}$. 
\end{lemma}
\begin{proof}
It is sufficient to prove that $\mathcal{K}(\Theta)=\mathsf{F}(V(0))$, 
since $\mathcal{H}(E)=E\mathcal{K}(\Theta)$. 
The proof below is essentially the same as the proof in \cite[Lemma 4.1]{Su20a}. 

If $\psi \in V(0)$, 
both  $\mathsf{F}\psi$ and $\mathsf{F}\mathsf{K}\psi$ belong to the Hardy space $H^2$ 
by definition \eqref{EQ_501} and $H^2=\mathsf{F}(L^2(0,\infty))$.   
On the other hand, we have 
$(\mathsf{F}\mathsf{K}\psi)(z)=\Theta(z)(\mathsf{F}\psi)^\sharp(z)$ 
by definition \eqref{EQ_501} again.  
This implies $(\mathsf{F}\psi)(z)=\Theta(z)(\mathsf{F}\mathsf{K}\psi)^\sharp(z)$, 
since $\Theta(z)\Theta^\sharp(z)=1$ by definition \eqref{EQ_104}.  
Therefore, $\mathsf{F}\psi$ belongs to $\mathcal{K}(\Theta)$ by \eqref{EQ_207}. 

Conversely, if $F \in \mathcal{K}(\Theta)$, there exists $f \in L^2(0,\infty)$ and $g \in L^2(-\infty,0)$ such that 
\begin{equation*}
F(z)=(\mathsf{F}f)(z) = \Theta(z) (\mathsf{F}g)(z). 
\end{equation*}
We have $(\mathsf{F}g)^\sharp(z) = \Theta(z)(\mathsf{F}f)^\sharp(z)$ by using $\Theta(z)\Theta^\sharp(z)=1$ again. 
Here $(\mathsf{F}g)^\sharp(z)=(\mathsf{F}\tilde{g})(z)$ for $\tilde{g}(x) = \overline{g(-x)} \in L^2(0,\infty)$, 
and $\Theta(z)(\mathsf{F}f)^\sharp(z)=(\mathsf{F}\mathsf{K}f)(z)$ as above. 
Hence $\mathsf{K}f$ belongs to $L^2(0,\infty)$, and thus $f \in V(0)$. 
\end{proof}
\begin{remark} \label{rmk_5_2}
By Lemma \ref{lem_5_1}, it follows that the RH would be false if $V(0) = {0}$, 
since $A(z)/(z-\gamma) = \xi(1/2 - iz)/(z-\gamma)$ belongs to $\mathcal{H}(E)$
for all $\gamma \in \Gamma$ under the assumption of the RH.
Therefore, it is an interesting problem to prove or disprove $V(0) \ne {0}$ unconditionally.
Since $V(0)$ is $\mathsf{K}$-invariant, if $V(0) \ne {0}$, then for any nonzero $f \in V(0)$,
the functions $(1 \pm \mathsf{K})f$ are eigenfunctions of $\mathsf{K}$
with eigenvalues $\pm 1$.
Hence, the problem reduces to determining whether the isometric involution $\mathsf{K}$
admits an eigenfunction in $L^2(0,\infty)$, which appears to be extremely difficult.
We therefore do not pursue this issue further in the present article. 
\end{remark}

Let $\tau=\tau_\xi$ be the measure on $\R$ determined from the screw function $g=g_\xi$ by 
\eqref{EQ_203}. 
Then, we have $g(0)=0$, $b=0$, and 
\begin{equation} \label{EQ_503}
d\tau(\lambda) = \sum_{\gamma\in\Gamma} 
m_\gamma \delta(\lambda-\gamma) \,d\lambda, \quad \lambda \in \R, 
\end{equation}
since 
\begin{equation*} 
g(t)=\sum_{\gamma \in \Gamma} m_\gamma \, \frac{e^{i\gamma t}-1}{\gamma^2}
\end{equation*}
by \cite[Theorem 1.1\,(2)]{Su22}, where $\delta$ is the Dirac mass at $\lambda=0$, 
%
We understand that 
the Hilbert space 
$L^2(\tau)$ is the space of sequences 
$S=(S(\gamma))_{\gamma \in \Gamma}$ with 
\begin{equation} \label{EQ_504} 
\Vert S \Vert_{L^2(\tau)}^2 
= 
\sum_{\gamma \in \Gamma} m_\gamma \,|S(\gamma)|^2. 
\end{equation}
Then, we prove two isomorphisms for $L^2(\tau)$ 
necessary for the proof of Theorem \ref{thm_1_1}.

\begin{lemma} \label{lem_5_2}
Hilbert spaces $V(0)$ and $L^2(\tau)$ are isomorphic by the linear map 
\[
V(0) \ni \psi ~\mapsto~ S_\psi := \Bigl( \widehat{\psi}(\gamma) \Bigr)_{\gamma \in \Gamma} 
\in L^2(\tau)
\]
with 
\begin{equation} \label{EQ_505}
2 \Vert \psi \Vert_{L^2(\R)}^2 = \Vert S_\psi \Vert_{L^2(\tau)}^2.
\end{equation}
\end{lemma} 
\begin{proof}
Let $\mu_\Theta$ be the measure on $\R$ determined from $\Theta=\Theta_\xi$ by \eqref{EQ_208}. 
Then, the linear map $\mathcal{K}(\Theta) \to L^2(\mu_\Theta)$ 
given by $\widehat{\psi} \mapsto S_\psi$ is an isometric isomorphism 
as reviewed in Section \ref{sec_2_4}. 
On the other hand, $L^2(\mu_\Theta)=L^2(\tau)$ 
with $\Vert S \Vert_{L^2(\mu_\Theta)}^2= \pi \Vert S \Vert_{L^2(\tau)}^2$
by \eqref{EQ_208}, \eqref{EQ_401}, and \eqref{EQ_503}. 
Therefore, by composing the maps $V(0) \to \mathcal{K}(\Theta)=\mathcal{F}(V(0))$ 
and $\mathcal{K}(\Theta) \to L^2(\mu_\Theta)$, we obtain the conclusion of the lemma, 
since $2 \pi \Vert \psi \Vert_{L^2(\R)}^2=\Vert \widehat{\psi} \Vert_{L^2(\R)}^2$. 
\end{proof}

\begin{lemma} \label{lem_5_3}
For $\displaystyle{\psi=\lim_{n\to\infty}\psi_n \in \mathcal{H}_W}$ with 
$\{\psi_n\}_{n \geq 1} \subset C_c^\infty(\R)$, 
we define $S_\psi \in L^2(\tau)$ by 
\[
S_\psi:=\lim_{n \to \infty} \Bigl( \widehat{\psi_n}(\gamma) \Bigr)_{\gamma \in \Gamma}
\quad \text{in} \quad L^2(\tau). 
\]
Then, it is well-defined 
and provides an isomorphism between $\mathcal{H}_W$ and $L^2(\tau)$ 
through the mapping
\[
\mathcal{H}_W \ni \psi ~\mapsto~ S_\psi \in L^2(\tau)
\]
with 
\begin{equation} \label{EQ_506}
\langle \psi, \psi \rangle_W = \Vert  S_\psi \Vert_{L^2(\tau)}^2 . 
\end{equation}
\end{lemma} 
\begin{proof}
We consider $C_0^\infty(\R)=\{\phi \in C_c^\infty(\R)\,|\, \widehat{\phi}(0)=0\}$, 
since we obtain the same completion $\mathcal{H}(G_g)$ 
even if we start from this space instead of $C_0(\R)$. 
Then differentiation $\psi \mapsto \psi'$ gives a bijection from $C_c^\infty(\R)$ to $C_0^\infty(\R)$. 
The inverse map is $\phi \mapsto \int_{-\infty}^{x} \phi(y)\,dy$. 
The Weil hermitian form and the hermitian form $\langle \cdot, \cdot \rangle_{G_g}$ 
defined by \eqref{EQ_202} for the screw function $g$ 
are related as in \eqref{EQ_410}, which is written as
\begin{equation} \label{EQ_507}
\langle \phi, \phi \rangle_{G_g} = \langle \psi, \psi \rangle_W, 
\quad \psi(x) = \int_{-\infty}^{x}\phi(y)\,dy, \quad \psi \in C_c^\infty(\R).
\end{equation}
(Although not necessary for the proof, 
$\langle \phi, \phi \rangle_{G_g}$ and $\langle \psi, \psi \rangle_W$ 
are positive definite on $C_0^\infty(\R)$ and $C_c^\infty(\R)$, respectively, 
by \cite[Lemma 2.1]{Su22}.)
Relation \eqref{EQ_507} extends to the completed Hilbert spaces. 
Therefore, $\mathcal{H}_W$ is isometrically isomorphic 
to the Hilbert space $\mathcal{H}(G_g)$ by $\mathcal{H}(G_g) \to \mathcal{H}_W :$ 
$\,[\phi] \mapsto [\psi]$ 
with $\psi=\lim_{n\to\infty}\psi_n$ and $\psi_n(x)=\int_{-\infty}^{x}\phi_n(y)\,dy$ 
for $\phi=\lim_{n\to \infty}\phi_n$ ($\phi_n \in C_c^\infty(\R)$). 

We define $\mathcal{H}(G_g) \to L^2(\tau)$ as follows. 
For $[\phi] \in \mathcal{H}(G_g)$, 
we define $S_\phi=(S_\phi(\gamma))_{\gamma \in \Gamma}$ $\in L^2(\tau)$ 
by 
\[
\lim_{n \to \infty} \Bigl( \widehat{\phi_n}(\gamma)/\gamma \Bigr)_{\gamma \in \Gamma}
\quad \text{in} \quad L^2(\tau)
\]
using a sequence $(\phi_n)_n$ in $C_0^\infty(\R)$ satisfying $\phi=\lim_{n\to\infty} \phi_n$. 
Then, the map is well-defined and 
$\langle [\phi], [\phi] \rangle_{G_g}
= \langle \phi, \phi \rangle_{G_g}
=\Vert S_\phi \Vert_{L^2(\tau)}$ by \eqref{EQ_202}, \eqref{EQ_409}, and \eqref{EQ_504}. 
Therefore, it establishes the isomorphic isomorphism 
$\mathcal{H}(G_g) \to L^2(\tau):$ $[\phi] \mapsto S_\phi$ 
\cite[Sections 5.3 and 12.5]{KrLa14}. 
Using $\mathcal{H}(G_g) \to \mathcal{H}_W$ and noting 
$\widehat{\phi}(\lambda)/\lambda=i\widehat{\psi}(\lambda)$ 
for $\phi \in C_0^\infty(\R)$, 
we define $\mathcal{H}_W \to L^2(\tau)$ by $[\psi] \mapsto S_\psi$ with 
%
\[
S_\psi
=(S_\psi(\gamma))_{\gamma\in\Gamma} 
=
\lim_{n \to \infty} \Bigl( \widehat{\psi_n}(\gamma) \Bigr)_{\gamma \in \Gamma}
= 
\lim_{n \to \infty} \Bigl( -i\widehat{\phi_n}(\gamma)/\gamma \Bigr)_{\gamma \in \Gamma}
\quad \text{in} \quad L^2(\tau), 
\]
where $(\phi_n)_n$ is a sequence in $C_0^\infty(\R)$ 
such that  $\psi=\lim_{n\to\infty}\psi_n$ with $\phi_n=\psi_n^\prime$. 
Then, the map is well-defined and 
\[
\langle [\psi], [\psi] \rangle_{W}= \langle \psi, \psi \rangle_{W}
=\Vert S_\psi \Vert_{L^2(\tau)}=\Vert S_\phi \Vert_{L^2(\tau)}
= \langle \phi, \phi \rangle_{G_g} = \langle [\phi], [\phi] \rangle_{G_g}
\] 
holds, where $\phi=\lim_{n\to \infty} \phi_n$ 
and the second equality follows from \eqref{EQ_102} and \eqref{EQ_504}. 
Hence, it establishes an isometric isomorphism 
$\mathcal{H}_W \to L^2(\tau)$ by $[\psi] \mapsto S_\psi$. 
As a result, the mapping $\mathcal{H}_W \to L^2(\tau)$  
is directly defined by  
$S_\psi=\lim_{n \to \infty} \Bigl( \widehat{\psi_n}(\gamma) \Bigr)_{\gamma \in \Gamma}$ 
and $[\psi] \mapsto S_\psi$ for $\psi=\lim_{n\to\infty}\psi_n$
with the desired equality for norms. 
\end{proof}

\begin{theorem} \label{thm_5_1} 
Assume that the RH is true. 
Let $\mathcal{H}_W$, $\mathcal{H}(E)$, and $\mathcal{K}(\Theta)$ be as above. 
Let $V(t)$ be the spaces defined in \eqref{EQ_502}. 
Then the following hold: 
\begin{enumerate}
\item $
\Vert E \widehat{\psi} \Vert_{\mathcal{H}(E)}^2
= \Vert \widehat{\psi} \Vert_{L^2(\R)}^2
= 2\pi \Vert \psi \Vert_{L^2(\R)}^2
= \pi \langle \psi, \psi \rangle_W
$ for $\psi \in V(0)$. 
\item 
The map from $\mathcal{K}(\Theta)$ to $\mathcal{H}_W$ 
obtained by the composition of the inverse of 
\begin{equation} \label{EQ_508}
V(0)~\to~\mathcal{K}(\Theta):~\psi ~\mapsto~\widehat{\psi}(z),  
\quad  2\pi\Vert \psi \Vert_{L^2(\R)}^2 = \Vert \widehat{\psi} \Vert_{L^2(\R)}^2 
\end{equation}
and
\begin{equation} \label{EQ_509}
V(0) ~\to~\mathcal{H}_W: ~ \psi ~\mapsto~[\psi], 
\quad 2 \Vert \psi \Vert_{L^2(\R)}^2 = \langle [\psi], [\psi] \rangle_W = \langle \psi, \psi \rangle_W
\end{equation}
agrees with the isomorphism $F \mapsto \psi_F$ in Theorem \ref{thm_1_1}. 
In particular, \eqref{EQ_509} is an isometric isomorphism up to a constant multiple. 
\end{enumerate}
\end{theorem}
\begin{proof} (1) It suffices to show that the equality 
\begin{equation} \label{EQ_510}
\Vert \psi \Vert_{L^2(\R)}^2
= \frac{1}{2} \langle \psi, \psi \rangle_W
\end{equation}
holds, since 
$ \Vert E \widehat{\psi} \Vert_{\mathcal{H}(E)}
= \Vert \widehat{\psi} \Vert_{L^2(\R)}
$ by \eqref{EQ_204} and 
$
\Vert \widehat{\psi} \Vert_{L^2(\R)}^2
= 2\pi \Vert \psi \Vert_{L^2(\R)}^2$ 
by \eqref{EQ_101}. 
For each $\gamma \in \Gamma$, we define $\psi_\gamma \in L^2(\R)$ by
\begin{equation} \label{EQ_511}
F_\gamma=\widehat{\psi_\gamma}.  
\end{equation} 
Then each $\psi_\gamma$ belongs to $V(0)$, and 
$\{\psi_\gamma\}_{\gamma \in \Gamma}$ forms an orthogonal basis 
satisfying  
\[
2\pi \Vert \psi_\gamma \Vert_{L^2(\R)}^2
= \Vert \widehat{\psi_\gamma} \Vert_{\mathcal{K}(\Theta)}^2 
=  \Vert F_\gamma \Vert_{\mathcal{K}(\Theta)}^2 
= 1
\] 
by Proposition \ref{prop_401} and Lemma \ref{lem_5_1}, 
since the orthogonality of $F_\gamma$'s is preserved under the Fourier transform. 
For $\psi=\sum_{\gamma} c_\gamma \psi_\gamma \in V(0)$, we have 
\[
\Vert \psi \Vert_{L^2(\R)}^2 = \frac{1}{2\pi} \sum_{\gamma \in \Gamma} |c_\gamma|^2
\]
by the orthogonality and 
\[
\widehat{\psi}(\gamma) = \frac{1}{\sqrt{m_\gamma \pi}}\, c_\gamma 
\]
by applying \eqref{EQ_402} to $\widehat{\psi}=\sum_\gamma c_\gamma F_\gamma$. 
From these two and \eqref{EQ_102}, we get \eqref{EQ_510}. 
\medskip

\noindent
(2) It is clear that the composition of the inverse of \eqref{EQ_508} and \eqref{EQ_509} 
agrees with the map $F \mapsto \psi_F$ of Theorem \ref{thm_1_1} 
including the equality for norms, and 
we observed in the proof of  Lemma \ref{lem_5_2} that 
the map \eqref{EQ_508} is an isometric isomorphism up to the  multiple $\sqrt{2\pi}$. 
Therefore, it suffices to show that the map \eqref{EQ_509} gives an isometric isomorphism 
up to the multiple $\sqrt{2}$. 

For $\psi \in V(0)$, the function $S_\psi \in L^2(\tau)$ is defined and satisfies  
$2 \Vert \psi \Vert_{L^2(\R)}^2 = \Vert S_\psi \Vert_{L^2(\tau)}^2$
by Lemma \ref{lem_5_2}. 
Then there exists a sequence $(\psi_n^\ast)_n \subset C_c^\infty(\R)$ 
that converges to $\psi^\ast$ with respect to $\langle\cdot,\cdot\rangle_W$ 
and  $S_\psi=S_{\psi^\ast}$ by Lemma \ref{lem_5_3}. 
The later implies 
$
\langle \psi- \psi_n^\ast, \psi - \psi_n^\ast \rangle_W
=
\langle \psi^\ast - \psi_n^\ast, \psi^\ast - \psi_n^\ast \rangle_W 
$ $\to 0$ ($n \to \infty$). 
Therefore, $\psi=\psi^\ast$, and 
hence $V(0) \to \mathcal{H}_W$ 
is directly defined by $\psi \mapsto [\psi]$. 
Furthermore, we obtain  
$2 \Vert\psi\Vert_{L^2(\R)}^2 = \langle \psi, \psi \rangle_W$ from 
\eqref{EQ_505} and \eqref{EQ_506}. 
Hence, this map is precisely the one given in \eqref{EQ_509}. 
\end{proof}

The equality 
$\Vert \psi \Vert_{L^2(\R)}^2 = 2^{-1} \langle \psi, \psi \rangle_W$ 
in Theorem \ref{thm_5_1} (1) shows that 
the $L^2$-structure induced from $L^2(\R)$ and 
``arithmetic structure'' (or ``local structure'')  
arising from the geometric side of the Weil explicit formula \eqref{EQ_303} 
coincide on a dense subspace of $V(0)$ 
consisting of functions for which the Weil explicit formula holds.

\begin{theorem} \label{thm_5_2} 
Let $\mathcal{H}_0$ and $\mathcal{K}_0$ are Hilbert spaces 
defined unconditionally in Section \ref{section_3_3}. 
Assume that the RH is true. 
Then, 
$\mathcal{H}_0=\mathcal{H}_W$ and $\mathcal{K}_0=\mathcal{K}(\Theta)$, 
and the extended map 
$\widehat{\mathcal{P}_D}: \mathcal{H}_W \to \mathcal{K}(\Theta)$ 
provides the inverse of the map in Theorem \ref{thm_5_1} (2).  
In particular, $V(0)$ is the $L^2$-closure of $V^\circ(0)$ in Corollary \ref{cor_1_1}. 
\end{theorem}
\begin{proof} 
For $\psi \in C_c^\infty(\R)$, we have 
\[
\Vert \widehat{\mathcal{P}_{D\psi}} \Vert_{L^2(\R)}^2
= \pi \sum_{\gamma \in \Gamma}  m_\gamma \vert\widehat{\psi}(\gamma)\vert^2
= \pi \langle \psi, \psi \rangle_W
\]
by \eqref{EQ_102}, \eqref{EQ_308}, and Proposition \ref{prop_401}. 
Hence, $\mathcal{H}_0$ coincides with $\mathcal{H}_W$ by definition \eqref{EQ_309}. 
Formula \eqref{EQ_308} shows that 
the image $\widehat{\mathcal{P}_{D\psi}}$ 
is defined independent of the representatives of $[\psi]$ in $\mathcal{H}_W$. 
On the other hand, $\mathcal{K}_0$ is a subspace of $\mathcal{K}(\Theta)$ 
by Proposition \ref{prop_401} again. 

We denote $F=\widehat{\mathcal{P}_{D\psi}}$ for $[\psi] \in \mathcal{H}_W$ 
and set $\psi_F = \mathsf{F}^{-1}(F)$ as in Theorem \ref{thm_1_1}. 
Then, $F(\gamma)= \widehat{\psi}(\gamma)$ by \eqref{EQ_308} and \eqref{EQ_402}. 
Therefore, $\widehat{\psi_F}(\gamma)=\psi(\gamma)$ 
for all $\gamma \in \Gamma$,  
and hence $[\psi]=[\psi_F]$ in $\mathcal{H}_W$. 
On the other hand, $\widehat{\mathcal{P}_{D\psi_F}}(z)=F$ by \eqref{EQ_308}, 
since $\widehat{\psi_F}=F$ by definition and $F(\gamma)= \widehat{\psi}(\gamma)$. 
Hence, we obtain the desired conclusion. 
\end{proof}

The totally ordered structure of the subspaces of the de Branges space $\mathcal{H}(E)$ 
is described by $V(t)$ as follows. 
 
\begin{theorem} \label{thm_5_3} 
Assume that the RH is true. 
Then, $E \,\mathsf{F}(V(t))$ is a de Branges subspace of $\mathcal{H}(E)$ 
for every $t \geq 0$ and is isometrically isomorphic to $\mathcal{H}_W(t)$ 
up to a constant multiple by the map of Theorem \ref{thm_1_1}. 
\end{theorem}
\begin{proof} It is sufficient to prove the first half of the theorem, 
since the second half follows from Theorem \ref{thm_5_1} (2). 
We prove the claim for positive $t$ such that $V(t) \not=\{0\}$, 
since the case of $t=0$ was proved in Lemma \ref{lem_5_1} 
and the claim is trivial if $V(t)=\{0\}$. 
The following is essentially the same as the proof of \cite[Lemma 4.3]{Su20a}. 

We show that $\mathcal{H}:=E(z)\mathsf{F}(V(t))$ 
is a Hilbert space consisting of entire functions and 
satisfies the axiom of the de Branges spaces: 
\begin{enumerate}
\item[(dB1)] For each $z \in \C\setminus \R$ the point evaluation $\Phi \mapsto \Phi(z)$ is a continuous linear functional on $\mathcal{H}$. 
\item[(dB2)] If $\Phi \in \mathcal{H}$, $\Phi^\sharp$ belongs to $\mathcal{H}$ and $\Vert \Phi \Vert_{\mathcal{H}} = \Vert \Phi^\sharp \Vert_{\mathcal{H}}$. 
\item[(dB3)] If $w \in \C \setminus \R$, $\Phi \in \mathcal{H}$ and $\Phi(w)=0$, 
\begin{equation*}
\frac{z-\bar{w}}{z-w}\Phi(z) \in \mathcal{H} \quad \text{and} \quad 
\left\Vert \frac{z-\bar{w}}{z-w}\Phi(z) \right\Vert_{\mathcal{H}} = \Vert \Phi \Vert_{\mathcal{H}},
\end{equation*}
\end{enumerate}
where the Hilbert space structure is the one induced from $V(t)$ 
that is equivalent to $\langle F,G\rangle_{\mathcal H}
= \int_{\R} F(z)\overline{G(z)}|E(z)|^{-2}dz$ for $F, G \in \mathcal{H}$. 

Let $\Phi(z)=E(z)(\mathsf{F}f)(z) \in \mathcal{H}$ with $f \in V(t)$. 
First, we prove that $\mathcal{H}$ consists of entire functions. 
We see that $\Phi(z)$ is holomorphic in $\C_+$ by $f \in L^2(t,\infty)$. 
If we write $(\mathsf{J}_\sharp f)(x):=\overline{f(-x)}$, 
the commutative relation $\mathsf{J}\mathsf{F}=\mathsf{F}\mathsf{J}_\sharp$ holds. 
Therefore, using \eqref{EQ_501} and $\mathsf{K}^2=1$, 
we have 
$\Phi(z)= E(z)(\mathsf{F}f)(z)= E^\sharp(z)(\mathsf{F}\mathsf{J}_\sharp \mathsf{K}f)(z)$. 
This shows that $\Phi(z)$ is also holomorphic in $\C_-$. 
Furthermore, 
$\mathsf{J}_\sharp \mathsf{K}f \in L^2(-\infty,-t)$, 
because  the tempered distribution kernel  $k:=\mathsf{F}^{-1}\Theta$ of $\mathsf{K}$ 
has support in $[0, \infty)$ by \cite[Theorems 1.1 and 1.2]{QXYYY09}. 
On the real line, 
$\lim_{z \to x}(\mathsf{F}f)(z)=(\mathsf{F}f)(x)$ and 
$\lim_{z \to x}(\mathsf{F}\mathsf{J}_\sharp \mathsf{K}f)(z)
=\lim_{z \to x}(\mathsf{F}\mathsf{K}f)^\sharp(z)=\Theta^\sharp(x)(\mathsf{F}f)(x)$ 
for almost all $x \in \R$, 
where $z$ is allowed to tend to $x$ nontangentially 
from $\C_+$ and $\C_-$, respectively. 
Hence, $\Phi(z)$ is also holomorphic in a neighborhood of each point of $\R$. 
By the above, $\Phi(z)$ is an entire function. 

We confirm (dB1). 
For $z \in \C_+$, $\Phi \mapsto \Phi(z)=E(z)\int_{t}^{\infty}f(x)e^{izx}dx$ is a continuous linear form. 
On the other hand, for $z \in \C_-$, $\Phi \mapsto \Phi(z)=E^\sharp(z)
\int_{-\infty}^{-t} \overline{(\mathsf{K}f)(-x)}e^{izx}\,dx$ is a continuous linear functional. 

We confirm (dB2). 
We have  
$\Phi^\sharp(z)
= E(z)(\mathsf{F}\mathsf{K}f)(z)$. 
Since $\mathsf{K}f \in V(t)$, 
the function $\Phi^\sharp$ belongs to $\mathcal{H}$. 
Since $\mathsf{K}$ is isometric, the equality of norms in (dB2) holds. 

We confirm (dB3). 
The equality of norms in (dB3) is trivial by the definition of 
the norm of $\mathcal{H}$. 
From (dB2), it is sufficient to show only the case of $w \in \C_+$. 
Suppose that $\Phi(w)=0$ for some $w \in \C_+$. 
Then $(\mathsf{F}f)(w)=0$, since $E(z)$ has no zeros on $\C_+$. 
We put $f_w(x):=f(x) - i(w-\bar{w}) \int_{0}^{x-t} f(x-y) e^{-iwy} dy$.   
Then we easily find that $f_w \in L^2(t,\infty)$ 
and $(\mathsf{F}f_w)(z)=((z-\bar{w})/(z-w))(\mathsf{F}f)(z)$ 
for $z \in \C_+$. 
Hence we complete the proof if it is shown that 
$\mathsf{K}f_w$ has support in $[t,\infty)$, 
since $\mathsf{K}f_w \in L^2(\R)$ by $f_w \in L^2(t,\infty)$. 
We put 
$g_w(x):=(\mathsf{K}f)(x) - i(\bar{w}-w) \int_{0}^{x-t} (\mathsf{K}f)(x-y) e^{-i\bar{w}y} dy$. 
Then $g_w$ has support in $[t,\infty)$ by $\mathsf{K}f \in L^2(t,\infty)$ 
and
$(\mathsf{F}g_w)(z)=((z-w)/(z-\bar{w}))(\mathsf{F}\mathsf{K}f)(z)
=(\mathsf{F}\mathsf{K}f_w)(z)$ for $z \in \C_+$.
Hence $g_w =\mathsf{K}f_w$ and the proof is completed. 
\end{proof}

We expect that $V(t) \ne {0}$ for some $t>0$, or rather that $V(t) \ne {0}$ holds for all $t \ge 0$, but we do not discuss this in the present article, 
since it seems to be a nontrivial problem related to the eigenfunctions of $\mathsf{K}$, 
as mentioned in Remark~\ref{rmk_5_2}.

\subsection{A weaker variant of Corollary \ref{cor_1_1}} 

Since the space $V(0)$ can be constructed unconditionally 
as well as $V^\circ(0)$ in Corollary \ref{cor_1_1}, 
it can be used to state an equivalence condition for the RH. 
However, since the construction of $V(0)$ is simpler than that of $V^\circ(0)$, 
more conditions are required for the equivalence condition.

\begin{proposition} \label{prop_5_1} 
Let $V(0)=L^2(0,\infty)\cap\mathsf{K}L^2(0,\infty)$ be as in \eqref{EQ_502}. 
Then the RH is true if and only if the following two conditions hold:  
\begin{enumerate}
\item $\Vert \psi \Vert_{L^2(\R)}^2 = 2^{-1} \langle \psi, \psi \rangle_W$ 
for every $\psi \in V(0)$.  
\item For a given $\gamma \in \Gamma$ and any $\epsilon>0$, 
there exists $\psi \in V(0)$ such that 
\[
\widehat{\psi}(-\gamma)=1, \qquad 
|\widehat{\psi}(-\gamma')| \leq \frac{\epsilon}{|\gamma-\gamma'|^{1+\delta}} \quad 
\text{for every $\gamma' \in \Gamma\setminus\{\gamma\}$}
\]
for some $\delta>0$ independent of $\gamma$, $\epsilon$, and $\psi$. 
\end{enumerate}
\end{proposition}
\begin{proof}
Assuming the RH, (1) follows from Theorem \ref{thm_5_1}\,(1). 
Also, (2) holds, since $\psi_\gamma=\mathsf{F}^{-1}(F_\gamma)$ in $V(0)$ satisfies 
$\widehat{\psi_{\gamma}}(\gamma)\not=0$ and 
$\widehat{\psi_\gamma}(\gamma')=0$ for $\gamma' \in \Gamma\setminus\{\gamma\}$.  

Conversely, we assume that (1) and (2) are satisfied. 
Then, we show that a contradiction arises if the RH is false. 
We take a nonreal $\gamma_0 \in \Gamma$. 
For any $\epsilon>0$, there exists $\psi_1$, $\psi_2\in V(0)$ 
such that $\widehat{\psi_1}(-\gamma_0)=i$, 
$\widehat{\psi_2}(-\overline{\gamma_0})=-i$, 
$|\widehat{\psi_1}(-\gamma)| \leq \epsilon |\gamma_0-\gamma|^{-1-\delta}$ 
for every $\gamma \in \Gamma\setminus\{\gamma_0\}$, 
and  
$|\widehat{\psi_2}(-\gamma)| \leq \epsilon |\overline{\gamma_0}-\gamma|^{-1-\delta}$ 
for every $\gamma \in \Gamma\setminus\{\overline{\gamma_0}\}$ by (2). 
Then, for $\psi:=\psi_1+\psi_2 \,(\not=0)$, we have 
$
\langle \psi, \psi \rangle_W 
= \sum_{\gamma \in \Gamma} m_\gamma 
\widehat{\psi}(-\gamma)(\widehat{\psi})^\sharp(-\gamma)
= -m_{\gamma_0} + O(\epsilon)
$,  
since $\sum_{\gamma \in \Gamma}m_\gamma|\gamma|^{-1-\delta}<\infty$. 
Therefore, $\langle \psi, \psi \rangle_W$ is negative for a sufficiently small $\epsilon>0$, 
but it contradicts (1). Hence the RH holds. 
\end{proof}

\comment{
\subsection{A variant of Theorem \ref{prop_5_1}} 

From the proof of Theorem \ref{prop_5_1}, 
it is easily shown that the following variant of Theorem \ref{prop_5_1} holds. 

\begin{theorem} \label{thm_4} 
The RH is true if and only if there exists a subspace $V$ of $L^2(\R)$ 
such that 
the following two conditions hold:  
\begin{enumerate}
\item $\Vert \psi \Vert_{L^2(\R)}^2 = 2^{-1} \langle \psi, \psi \rangle_W$ 
for every $\psi \in V$.  
\item For a given $\gamma \in \Gamma$ and any $\epsilon>0$, 
there exists $\psi \in V$ such that 
\[
\widehat{\psi}(\gamma)=1, \qquad 
|\widehat{\psi}(\gamma')| \leq \frac{\epsilon}{|\gamma-\gamma'|^{1+\delta}} \quad 
\text{for every $\gamma' \in \Gamma\setminus\{\gamma\}$}
\]
for some $\delta>0$ independent of $\gamma$, $\epsilon$, and $\psi$. 
\end{enumerate}
\end{theorem}
}

\section{Hilbert--P{\'o}lya space} \label{section_6}

One of attractive strategies for proving the RH 
is the construction of a Hilbert--P{\'o}lya space, 
which is a pair of a Hilbert space and a self-adjoint operator acting on it 
such that all nontrivial zeros of the Riemann zeta-function 
are eigenvalues of the self-adjoint operator. 
In this section, 
we state that $\mathcal{H}_W$ is one of Hilbert--P{\'o}lya spaces {\it under the RH}. 
Note that $\mathcal{H}_W$ is unconditionally defined as $\mathcal{H}_0$ 
by Theorem \ref{thm_5_2}. 
\smallskip

We assume the RH and denote $E=E_\xi$ as in Section \ref{section_5}. 
In this case, the domain $\mathfrak{D}(\mathsf{M})$ of 
the multiplication operator $\mathsf{M}$ on $\mathcal{H}(E)$
is dense in $\mathcal{H}(E)$, 
because $S_\theta(z)$ does not belongs to $\mathcal{H}(E)$ for all $\theta \in [0,\pi)$ 
by the estimate $|S_\theta(iy)/E(iy)| \gg (\log y)^{-1}$ ($y \to +\infty$) 
obtained by the Stirling formula for the gamma-function 
and \cite[Proposition 2.1]{Re02}. 
Using $\mathsf{M}$, 
we define the operator $\mathsf{A}:= \mathsf{F}^{-1}\mathsf{M} \mathsf{F}$ on $V(0)$ 
with the domain $\mathfrak{D}(\mathsf{A})=\mathsf{F}^{-1}(\mathfrak{D}(\mathsf{M}))$. 
If $\psi \in V(0)$ is differentiable and $\psi'$ also belongs to $V(0)$, 
then $\mathsf{A}\,\psi=i\psi'$. 
Further, we define the operator $\mathsf{A}_W$ on $\mathcal{H}_W$ as follows. 
\smallskip

By Theorem \ref{thm_5_2}, 
the inverse of \eqref{EQ_509} from $\mathcal{H}_W$ to $V(0)$ 
is given by $[\psi] \mapsto \mathsf{F}^{-1}\widehat{\mathcal{P}_{D\psi}}$. 
Further, if we choose the representative of $\psi$ from $V(0)$, 
it is possible and uniquely determined by Theorem \ref{thm_5_1} (2), 
and therefore $\psi=\mathsf{F}^{-1}\widehat{\mathcal{P}_{D\psi}}$.  
By choosing representatives in this way, 
we define $\mathsf{A}_W$ on $\mathcal{H}_W$ 
by $\mathsf{A}_W[\psi]=[\mathsf{A}\psi]$.   
By the same procedure as above, 
the family of self-adjoint extensions $\mathsf{M}_{\theta}$ of $\mathsf{M}$ 
determines the corresponding families 
of self-adjoint extensions of $\mathsf{A}$ and $\mathsf{A}_W$ (see \eqref{EQ_205} and \eqref{EQ_206}). 
By this correspondence, 
the orthogonal basis $\{[\psi_\gamma]\}_{\gamma \in \Gamma}$ of $\mathcal{H}_W$ 
consists of eigenvectors $[\psi_\gamma]$ of $\mathsf{A}_{W,\pi/2}$ 
with eigenvalues $\gamma \in \Gamma$, 
since $\{E F_\gamma\}_{\gamma \in \Gamma}$ with \eqref{EQ_305} 
is an orthogonal basis of $\mathcal{H}(E)$ 
consists of eigenfunctions of $\mathsf{M}_{\pi/2}$ 
with eigenvalues $\Gamma$ (see Seciton \ref{sec_2_3}). 
Therefore, the pair $(\mathcal{H}_W,\,\mathsf{A}_{W,\pi/2})$ 
is a Hilbert--P{\'o}lya space.

It is important to note that the multiplicity of $\gamma \in \Gamma$ as an eigenvalue of 
$\mathsf{A}_{W,\pi/2}$ (and $\mathsf{M}_{\pi/2}$) is one. 
In other words, the multiplicity of $\gamma \in \Gamma$ 
as a zero of $\xi(1/2-iz)$ is not reflected in the multiplicity of 
$\mathsf{A}_{W,\pi/2}$ (and $\mathsf{M}_{\pi/2}$). 
In particular, 
it shows the explicit difference between the de Branges space $\mathcal{H}(E_\xi)$ 
and the de Branges space $\mathcal{B}_\lambda^{S}$ in \cite[Section 4.8]{CCM24}. 
\medskip

In the above discussion, we assumed the RH, 
but \eqref{EQ_205} and \eqref{EQ_206} allow us to define the operator 
$\mathsf{M}_\theta$ without the RH. 
However, its properties as an operator become unclear.

\section{Special values of the screw line $\mathfrak{S}_t(z)$} \label{Section_4} 

The screw line $\mathfrak{S}_t(z)$ 
has the following unconditional relations with the screw function $g(t)$. 
It is interesting that they are not a special case of equations 
obtained from the general theory of screw functions.

\begin{theorem} 
Let $g_\xi(t)$ and $\mathfrak{P}_t(z)$ 
be functions of \eqref{EQ_403} and \eqref{EQ_106}, respectively.  
Then the following equations hold independently of the truth of the RH: 
\begin{equation} \label{EQ_701}
\mathfrak{P}_t(0) = -g_\xi(t), \\
\end{equation}
\begin{equation} \label{EQ_702}
\lim_{y\to + \infty} \left[ y\, \mathfrak{B}_t(-iy) -\frac{1}{2}
\frac{\Gamma'}{\Gamma}\left(\frac{1}{4}+\frac{y}{2} \right) 
+\frac{1}{2} \log \pi
\right] = -g_\xi^\prime(t),  
\end{equation}
where we assume $t \not = \log n$ for any $n \in \N$ in \eqref{EQ_702}.  
\end{theorem}

\begin{proof} 
Equality \eqref{EQ_701} follows 
from \eqref{EQ_302}, Proposition \ref{prop_301}, and 
\cite[Theorem 1.1 (2)]{Su22}, 
but it follows directly from \eqref{EQ_403} and \eqref{EQ_106} as follows.  
By $\Phi(z,s,a)=\sum_{n=0}^{\infty}z^n(n+a)^{-s}$ 
and \eqref{EQ_208}, 
\[
\lim_{z\to 0}
\frac{1}{iz}
\Bigl[ \Phi(e^{-2t},1,\tfrac{1}{2}(\tfrac{1}{2}-iz)) - \Phi(e^{-2t},1,1/4) \Bigr] 
= -\frac{1}{2} \Phi(e^{-2t},2,1/4), 
\]
\[
\lim_{z\to 0}
\frac{1}{iz}\left[
\frac{\Gamma'}{\Gamma}\left(\frac{1}{4}-\frac{iz}{2}\right)
-
\frac{\Gamma'}{\Gamma}\left(\frac{1}{4}\right)
\right]
= \frac{1}{2}\psi_1\left(\frac{1}{4}\right), 
\]
where $\psi_1(z)$ is the polygamma function of order one.  
The expansion 
$\psi_1(w)=\sum_{n=0}^{\infty}(w+n)^{-2}$ 
gives $\psi_1(1/4)=\Phi(1,2,1/4)$. 
Taking  $s=1/2$ in the logarithmic derivative of $\xi(s)=\xi(1-s)$ 
and using
\[
\frac{\Gamma'}{\Gamma}\left(\frac{1}{4}\right) 
= -\gamma_0 - 3\log 2 -\frac{\pi}{2}, 
\]
we have 
\[
\frac{\zeta'}{\zeta}\left(\frac{1}{2}\right)
= \frac{1}{2} \left(\gamma_0+3\log 2+\log \pi + \frac{\pi}{2} \right). 
\]
Hence, by taking the limit $z \to 0$ in \eqref{EQ_106}, 
we obtain the minus of \eqref{EQ_403}. 

To show \eqref{EQ_702}, we multiply \eqref{EQ_106} 
by $y$ and substitute $-iy$ for $z$: 
\[
\aligned
y \, \mathfrak{P}_t(-iy)
& := \frac{4y(e^{t/2}-1)}{1+2y} + \frac{4y(e^{-t/2}-1)}{1-2y} \\ 
& \quad 
+ (e^{-yt}-1)\frac{\zeta'}{\zeta}\left( \frac{1}{2}-y \right)
+ \sum_{n \leq e^t} \frac{\Lambda(n)}{\sqrt{n}} (e^{-y(t-\log n)}-1)
\\
& \quad 
+ \frac{1}{2}
\left[
\frac{\Gamma'}{\Gamma}\left(\frac{1}{4}\right)
-
\frac{\Gamma'}{\Gamma}\left(\frac{1}{4}-\frac{y}{2}\right)
\right] \\
& \quad 
+ \frac{1}{2} 
e^{-t/2} \left[
\Phi(e^{-2t},1,1/4)
-
\Phi(e^{-2t},1,\tfrac{1}{2}(\tfrac{1}{2}-y)) 
\right]
\endaligned 
\]
Therefore, for positive $t>0$, 
\[
\aligned 
\,&\lim_{y\to + \infty}  \left[ y\, \mathfrak{B}_t(-iy) -\frac{1}{2}
\frac{\Gamma'}{\Gamma}\left(\frac{1}{4}+\frac{y}{2} \right) 
+\frac{1}{2} \log \pi
\right] \\
& = 2(e^{t/2}-e^{-t/2})  -  \sum_{n \leq e^{t}} \frac{\Lambda(n)}{\sqrt{n}} 
+ \frac{1}{2}\left[
\frac{\Gamma'}{\Gamma}\left(\frac{1}{4}\right)-\log\pi \right]+\frac{1}{2}e^{-t/2}\Phi(e^{-2t},1,1/4)
\endaligned 
\]
by using the logarithmic derivative of $\xi(s)=\xi(1-s)$ at $s=1/2-y$. 
The right-hand side equals $-g_\xi^\prime(t)$ if $t\not=\log n$ 
by \eqref{EQ_403}, and 
$(d/dt)(e^{-t/2}\Phi(e^{-2t},2,1/4))=-2e^{-t/2}\Phi(e^{-2t},2,1/4)$ 
follows from $\Phi(z,s,a)=\sum_{n=0}^{\infty}z^n(n+a)^{-s}$. 
\end{proof}

\medskip

\noindent
{\bf Acknowledgment} 
\medskip

\noindent
This work was supported by JSPS KAKENHI Grant Number 
JP23K03050. 
This work was also supported by the Research Institute for Mathematical Sciences, 
an International Joint Usage/Research Center located in Kyoto University.

%

%
\bigskip 

\noindent
Masatoshi Suzuki,\\[5pt]
Department of Mathematics, \\
School of Science, \\ 
Institute of Science Tokyo \\
2-12-1 Ookayama, Meguro-ku, \\
Tokyo 152-8551, Japan  \\[2pt]
Email: {\tt msuzuki@math.sci.isct.ac.jp}


\begin{thebibliography}{99}
%
\bibitem{BS08}
N. A. Baas, C. F. Skau, 
\newblock{The lord of the numbers, Atle Selberg. On his life and mathematics}, 
\newblock{\it Bull. Amer. Math. Soc. (N.S.)} 
\newblock{{\bf 45} (2008), no. 4, 617--649}. 
%
\bibitem{BoLa99}
E. Bombieri, J. C. Lagarias, 
\newblock{Complements to Li's criterion for the Riemann hypothesis}, 
\newblock{\it J. Number Theory} 
\newblock{{\bf 77} (1999), no. 2, 274--287}. 
%
\bibitem{Bo01}
E. Bombieri,
\newblock{Remarks on Weil's quadratic functional in the theory of prime numbers. I}, 
\newblock{\it Atti Accad. Naz. Lincei Cl. Sci. Fis. Mat. Natur. Rend. Lincei (9) Mat. Appl.} 
\newblock{{\bf 11} (2000), no. 3, 183--233 (2001)}. 
%
\bibitem{CCM24}
A. Connes, C. Consani, H. Moscovici, 
\newblock{Zeta zeros and prolate wave operators}, 
\newblock{\it  Annals of Functional Analysis} 
\newblock{{\bf 15} (2024), article number 87}. 
%
\bibitem{dB68}
L. de Branges,
\newblock{Hilbert spaces of entire functions}, 
\newblock{\it Prentice-Hall, Inc., Englewood Cliffs, N.J.}   
\newblock{1968}.
%
\bibitem{KW99}
M. Kaltenb{\"a}ck, H. Woracek, 
\newblock{Pontryagin spaces of entire functions. {I}},
\newblock{\it Integral Equations Operator Theory }  
\newblock{{\bf 33} (1999), no. 1, 34--97}. 
%
\bibitem{KrLa14}
M. G. Kre\u{\i}n, H. Langer, 
\newblock{Continuation of hermitian positive definite functions and related questions}, 
\newblock{\it Integral Equations Operator Theory}  
\newblock{{\bf 78} (2014), no. 1, 1--69}. 
%
\bibitem{La06}
J. C. Lagarias,
\newblock{Hilbert spaces of entire functions and Dirichlet $L$-functions}, 
\newblock{\it Frontiers in number theory, physics, and geometry. I},   
\newblock{365--377}, 
\newblock{\it Springer, Berlin}, 
\newblock{2006}.
%
\bibitem{MaPo05}
N. Makarov, A. Poltoratski, 
\newblock{Meromorphic inner functions, Toeplitz kernels and the uncertainty principle}, 
\newblock{\it Perspectives in analysis},   
\newblock{185--252}, 
\newblock{Math. Phys. Stud., 27}, 
\newblock{\it Springer, Berlin}, 
\newblock{2005}.
%
\bibitem{QXYYY09}
T. Qian, Y. Xu, D. Yan, L. Yan, B. Yu, 
\newblock{Fourier spectrum characterization of Hardy spaces and applications}, 
\newblock{\it Proc. Amer. Math. Soc.}   
\newblock{{\bf 137} (2009), no. 3, 971--980}. 
%
\bibitem{Re02}
C. Remling, 
\newblock{Schr\"{o}dinger operators and de {B}ranges spaces}, 
\newblock{\it J. Funct. Anal.}  
\newblock{{\bf 196} (2002), no. 2, 323--394}. 
%
\bibitem{SiTo15}
L. O. Silva, J. H. Toloza, 
\newblock{de Branges spaces and Kre\u{\i}n's theory of entire operators}, 
\newblock{\it Operator Theory}, 
\newblock{D. Alpay (eds.)},
\newblock{\it Springer, Basel}, 
\newblock{2015}, 
\newblock{pp. 549--580}. 
%
\bibitem{Su20a} 
M. Suzuki,
\newblock{An inverse problem for a class of canonical systems having Hamiltonians of determinant one}, 
\newblock{\it J. Funct. Anal.},  
\newblock{{\bf 279} (2020), no. 12, 108699}. 
%
\bibitem{Su22} 
M. Suzuki,
\newblock{Aspects of the screw function corresponding to the Riemann zeta function}, 
\newblock{\it J. Lond. Math. Soc.}  
\newblock{{\bf 108} (2023), no.4, 1448-1487}. 
%
\bibitem{Su23}
M. Suzuki,
\newblock{Li coefficients as norms of functions in a model space}, 
\newblock{\it J. Number Theory}  
\newblock{{\bf 252} (2023), 177--194}. 
%
\bibitem{Tit86}
E. C. Titchmarsh,
\newblock{The theory of the Riemann zeta-function, Second edition, 
Edited and with a preface by D. R. Heath-Brown
}, 
\newblock{\it The Clarendon Press, Oxford University Press, New York}, 
\newblock{1986}. 
%
\bibitem{We52}
A. Weil,
\newblock{Sur les ``formules explicites'' de la th\'{e}orie des nombres
              premiers}, 
\newblock{\it Comm. S\'{e}m. Math. Univ. Lund [Medd. Lunds Univ. Mat. Sem.]},  
\newblock{{\bf 1952} (1952), Tome Suppl\'{e}mentaire, 252--265}. 
%
\bibitem{Wo15}
H. Woracek, 
\newblock{De Branges spaces and growth aspects}, 
\newblock{\it Operator Theory}, 
\newblock{D. Alpay (eds.)},
\newblock{\it Springer, Basel}, 
\newblock{2015}, 
\newblock{pp. 489--523}. 
%
\bibitem{Yo92}
H. Yoshida,
\newblock{On Hermitian forms attached to zeta functions}, 
\newblock{\it Zeta functions in geometry (Tokyo, 1990)}, 
\newblock{281--325, Adv. Stud. Pure Math., 21}, 
\newblock{\it Kinokuniya, Tokyo}, 
\newblock{1992}.
%
\end{thebibliography}
\end{document}